\documentclass[aos]{imsart}
\usepackage[sort]{natbib}
\usepackage{fancyhdr}
\usepackage{mathtools}
\usepackage{bbm}
\usepackage{color}
\usepackage{listings}
\usepackage{array}
\usepackage{mathrsfs}
\newcolumntype{P}[1]{>{\centering\arraybackslash}p{#1}}

\usepackage{multirow}
\usepackage{xr}

\newcommand{\ssymbol}[1]{^{\@fnsymbol{#1}}}
\makeatother
\usepackage[utf8]{inputenc}
\usepackage[T1]{fontenc}
\usepackage{lmodern}
\usepackage{graphicx}
\usepackage{floatrow}
\usepackage{amsmath,amssymb,amsthm}
\usepackage[english]{babel}
\newcommand{\RN}[1]{%
	\textup{\uppercase\expandafter{\romannumeral#1}}%
}

\graphicspath{{plots/}}

\newcommand{\bbone}{\mathbbm{1}}

\def\mH{\mathcal{H}}
\def\mI{\mathbb{I}}
\def\mB{\mathbb{B}}
\def\mU{\mathcal{U}}
\def\Re{\mathbb{R}}

\newcommand{\overbar}[1]{\mkern 1.5mu\overline{\mkern-1.5mu#1\mkern-1.5mu}\mkern 1.5mu}
\def\bd{\boldsymbol{d}}
\def\bc{\boldsymbol{c}}
\def\bg{\boldsymbol{g}}
\def\bK{\overbar{K}}
\def\ka{\kappa}
\def\eps{\epsilon}
\def\ti{\tilde}
\def\bXi{\boldsymbol{\Xi}}

\def\bdel{\boldsymbol{\delta}}
\def\utau{\underline{\tau}}
\def\otau{\overline{\tau}}
\def\convd{\stackrel{d}{\longrightarrow} }

\newcommand\ca[1]{{\cal{#1}}}
\newcommand\lo[1]{_{\nano #1}}
\def\nano{\scriptscriptstyle}
\newcommand\hi[1]{^{\nano #1}}
\newcommand\inner[2]{\langle #1, #2 \rangle}
\newcommand{\trans}{^{\mbox{\tiny {\sf T}}}}
\newcommand{\norm}[1]{\left\lVert#1\right\rVert}

\theoremstyle{plain}
\newtheorem{proposition}{{\bf Proposition}}
\numberwithin{proposition}{section}
\newtheorem{lemma}{{\bf Lemma}}
\numberwithin{lemma}{section}
\newtheorem{corollary}{{\bf Corollary}}
\numberwithin{corollary}{section}
\newtheorem{theorem}{{\bf Theorem}}
\numberwithin{theorem}{section}

\theoremstyle{definition}

\newtheorem{remark}{{\bf Remark}}
\numberwithin{remark}{section}
\newtheorem{assump}{Assumption}

\DeclareMathOperator*{\argmin}{arg\,min}

\numberwithin{equation}{section}

\def\Var{{\rm Var}\,}
\def\E{{\rm E}\,}

\begin{document}
	
	\begin{frontmatter}
		\title{Statistical Inference for Functional Linear Quantile Regression}
		\runtitle{FLQR}
		
		\begin{aug}
			\author[A]{\fnms{Peijun} \snm{Sang}\ead[label=e1]{psang@uwaterloo.ca}}, 
			\author[B]{\fnms{Zuofeng} \snm{Shang}\ead[label=e2]{zshang@njit.edu}} and
			\author[C]{\fnms{Pang} \snm{Du}\ead[label=e3]{pangdu@vt.edu}}
			\address[A]{Department of Statistics and Actuarial Science,
				University of Waterloo, 
				\printead{e1}}
			
			\address[B]{Department of Mathematical Sciences, New Jersey Institute of Technology,
				\printead{e2}}
			
			\address[C]{Department of Statistics,
				Virginia Polytechnic Institute and State University,  
				\printead{e3}}
			
		\end{aug}

		\begin{abstract}
			We propose inferential tools for functional linear quantile regression  where the conditional quantile of a scalar response is assumed to be a linear functional of a functional covariate. In contrast to conventional approaches, we employ kernel convolution to smooth the original loss function. The coefficient function is estimated under a reproducing kernel Hilbert space framework. A gradient descent algorithm is designed to minimize the smoothed loss function with a roughness penalty.
			With the aid of the Banach fixed-point theorem, we show the existence and uniqueness of our proposed estimator as the minimizer of the regularized loss function in an appropriate Hilbert space. 
			Furthermore, we establish the convergence rate as well as the
			weak convergence of our estimator. As far as we know, this is the first weak convergence result for a functional quantile regression model. 
			Pointwise confidence intervals and a simultaneous confidence band for the true coefficient function are then developed based on these theoretical properties. 
			Numerical studies including both simulations and a data application are conducted to investigate the performance of our estimator and inference tools in finite sample. 
		\end{abstract}
		
		\begin{keyword}[class=MSC2020]
			\kwd[Primary ]{00X00}
			\kwd{00X00}
			\kwd[; secondary ]{00X00}
		\end{keyword}
		
		\begin{keyword}
			reproducing kernel Hilbert space, Functional Bahadur representation, Pointwise confidence intervals, Weak convergence
		\end{keyword}
		
	\end{frontmatter}
	
	\section{Introduction}
	Functional data analysis (FDA), extending the traditional data domain to curves, surfaces, and objects, has attracted a lot of attention in the past decades. Monographs, such as \cite{ramsay2005}, \cite{ferraty2006}, \cite{horvath2012}, \cite{hsing2015} and \cite{kokoszka2017}, exhibit a comprehensive review of the past and current topics in FDA. Among all the topics, functional linear models focusing on the mean response or the mean response function are arguably the most studied. Theoretical investigation on the asymptotic consistencies in coefficient function estimation and mean prediction has been mainly considered under two frameworks: the functional principal component analysis (FPCA) and the reproducing kernel Hilbert space (RKHS) frameworks. Examples for FPCA include \cite{yao2005b}, \cite{cai2006}, and \cite{hall2007}. Examples for the RKHS approach include \cite{yuan2010}, \cite{cai2012}, and \cite{sun2018optimal}. In this paper, we shall study, under the RKHS framework, the functional linear quantile regression (FLQR) model where the modeling of the (conditional) quantile function of a scalar response against a functional covariate is of interest.
	
	Comparing with various functional linear models, functional quantile regression models have been studied much less in the literature with some exceptions given below. \cite{cardot2005} considered the convergence of a penalized B-spline coefficient function estimator under a functional quantile regression model. \cite{ferraty2005} and \cite{chen2012} estimated the conditional quantile function through inverting an estimate of the conditional cumulative distribution of the response given the functional covariate. The difference was that the former used kernel estimation and the latter used a generalized functional linear model via FPCA. 
	\cite{kato2012} considered a more thorough and direct approach. The conditional quantile function over a range of quantile indices was modeled through a function-on-function regression model against a functional covariate. The FPCA was used to estimate the coefficient function and the resulting estimator of the conditional quantile function was further monotonized to satisfy the monotonicity constraint over the quantile index. A minimax optimal rate of convergence was established for the estimators. \cite{yu2016} considered a model similar to \cite{cardot2005} but the basis functions used for estimating the coefficient function were generated from a so-called partial quantile regression technique that resembles the partial least squares. They also extended the approach to the compositive quantile regression scenario \citep{zou2008}. \cite{yao2017} considered quantile regression models with a functional covariate and a high dimensional predictor variable where the focus was on variable selection. \cite{zhu2019} considered a similar high dimensional model with multiple functional covariates involved and double penalization applied to select both functional and scalar variables. However, one common drawback in these papers is the lack of any rigorously derived weak convergence results and inference tools.
	
	Such kind of inference theory has been studied for functional regression models with a mean response.
	To test on the nullity of the slope function in a functional linear model, \cite{cardot2003} proposed two test statistics based on the norm of the empirical cross-covariance operator of the response and the functional predictor with a pre-selected number of FPCs of the predictor process.
	\cite{muller2005} derived the asymptotic distribution for the coefficient function estimate under a generalized functional linear regression model fitted by the FPCA, where the truncation point of the Karhunen-Lo\`{e}ve expansion is assumed to increase at a certain rate of the sample size and randomness of FPCs is ignored. 
	\cite{cardot2007} developed a central limit theorem for the FPCA approach under a functional linear regression model where the predictor function can reside in a more general Hilbert space. 
	\cite{ferraty2007} and \cite{ferraty2010} considered nonparametric functional regression models based on kernel estimation. They derived point-wise normality results for the regression function and developed bootstrap procedures for the construction of point-wise confidence intervals.   A similar bootstrap procedure was proposed in \cite{gonz2011} for functional linear models. \cite{zhang2007} considered statistical inference for functional linear models where the functional data are reconstructed by local polynomial kernel estimation and proposed a global Wald type of test statistic on the effect of functional covariates. To test the nullity of the slope function in a functional linear model, \cite{hilgert2013} studied a Fisher-type nonadaptive test statistic corresponding to the projection of the response on the leading FPCs of the predictor process where the randomness of the FPCs are incorporated. Employing a sequential approximation by a series of functional principal components regression models, \cite{lei2014} proposed a global test procedure on the slope function in a functional linear model. Under the framework of RKHS, \cite{shang2015} developed a system of statistical inference tools for generalized functional linear models, which include confidence intervals/bands, prediction intervals, functional contrast tests, and global tests on slope functions. Extensions of their approach to a functional Cox model and a function-on-function regression model were considered respectively in \cite{hao2020} and \cite{dette2021}. \cite{cuesta2019} constructed goodness-of-fit tests for the functional linear model with a scalar response based on marked empirical processes indexed by a randomly projected functional covariate. 

	Statistical inference for functional quantile regression models has been rarely studied, partly due to the fact that the quantile loss function is non-differentiable. To our best knowledge, only one informal treatment was attempted in \cite{li2016inference} where an adjusted Wald test was proposed to test on whether a common coefficient function is shared across models for a given set of multiple quantile indices. However, a key assumption in their work is that the coefficient function can be represented by the expansion of a {\it fixed number} of functional principle components. This strong assumption essentially reduces the problem to the traditional quantile linear regression model. Then all the traditional theoretical tools can be directly applied to yield the root-$n$ convergence rate of the coefficient estimate and an asymptotic chi-square distribution for the test statistic. The problem considered here is much harder and the treatment here is more rigorous than those in \cite{li2016inference} from the following aspects. First, we don't have a fixed dimension assumption for the coefficient function which we assume belongs to an infinite dimensional Sobolev space. It is well known that statistical inference on an infinite dimensional parameter space is a much harder problem. Second, we consider the model over a continuum range of quantile indices rather than the composite quantile regression scenario that only entertains a small number of quantile indices. Such a functional quantile regression model has only been studied by \cite{kato2012} with a minimax estimation consistency result.
	
	Therefore, this paper aims to provide a systematic and rigorous study on the inference problem for a FLQR model with a scalar response and a functional covariate over a range of quantile indices. Recognizing the non-differentiability of the quantile loss function, we first employ a kernel density smoothing technique introduced in \cite{fernandes2019} to consider a smoothed loss function. With the addition of a roughness penalty to the loss function, our coefficient function estimate is defined as the minimizer of this penalized and smoothed loss function. The representer theorem guarantees that this minimizer lies on a finite dimensional subspace (note that the finite dimension is not fixed but increases with $n$), although the Sobolev space of the original optimization problem is of infinite dimensions. Based on this we develop an efficient gradient descent algorithm to compute the coefficient function estimate. For the theoretical properties, we first show that the minimizer of the penalized and smoothed objective function provides a good approximation to the minimizer of the penalized quantile loss objective function with a negligible error.  Then a functional Bahadur representation is derived for the smoothed FLQR model estimate. Based on the representation, the weak convergence of the coefficient function estimate to a Gaussian process is established. We then derive the pointwise confidence interval and the simultaneous confidence band for the coefficient function, as well as the confidence interval for the conditional quantile. The techniques we use in this paper are very different from those in the earlier work on penalized functional regression models, such as \cite{yuan2010} and \cite{cai2012}, where the alignment of the covariance function and the reproducing kernel plays an important role. Useful as they were in deriving the optimal estimation and prediction consistencies, the covariance alignment argument can be difficult to extend to the inference scenario. Here, we rely on the Banach fixed-point theorem to derive the functional Bahadur representation. The work here is also significantly different from that in \cite{shang2015} who studied the functional Bahadur representation for a generalized functional linear model. Comparing with the smooth and convex objective function in their model setting, the quantile loss function here is non-differentiable. Although the smoothed objective function trick makes the computation more tractable, considerable effort is needed to fill in the gap between the smoothed approximation and the original objective settings. 
	
	In summary, our work makes the following contributions to FLQR models. First, we provide a smoothed approximation to the original FLQR model with a negligible approximation error. Second, we develop an efficient estimation procedure whose resulting coefficient function estimator is not only theoretically consistent but also delivers better accuracy than the existing methods in empirical experiments. Lastly and more importantly, we derive the functional Bahadur representation for the coefficient function estimator, based on which the inference tools, such as simultaneous confidence bands for the coefficient function and confidence intervals for the conditional quantile, are constructed. As far as we know, these inference tools are the first ones for a functional quantile regression model.
	
	The rest of the paper is structured as follows. In Section \ref{sec:model} we introduce the smoothed version of FLQR and the gradient descent method to fit the smoothed FLQR with a roughness penalty. Theoretical properties such as consistency, convergence rate and weak convergence of our proposed estimator are discussed in Section \ref{sec:Bahadur}. In Sections \ref{sec:simulation} and \ref{sec:real} we investigate finite sample performance of the proposed estimator as well as inference tools through empirical studies. Section \ref{sec:conclusion} concludes the article. All technical proofs are delegated to the supplementary material. 
	
	\section{Model and Estimation}
	\label{sec:model}
	Let $X$ denote a random function defined on $\mI$, a compact subset of $\Re$ and $Y$ be a scalar random variable taking values in $\ca Y$. Without loss of generality, we assume that $\mI = [0, 1]$. 
	Let $Q_{Y|X}(\cdot | X)$ denote the conditional quantile of $Y$ given $X$. Throughout the paper the conditional quantile of $Y$ w.r.t. $X = X(\cdot)$ can be understood as a function of a collection
	of random variables $\{X(t): 0 \leq t \leq 1\}$. 
	Let $\ca U$ be a given subset of (0, 1) that is away from 0 and 1, that is , for some small $c_0 \in (0, 1/2),~ \mU \subset [c_0, 1 - c_0]$. 
	We consider the following functional linear quantile regression (FLQR) model: 
	\begin{equation}
		Q_{Y|X}(\tau | X) = \alpha_0(\tau) + \int_0^1 \beta_0(t, \tau)X(t) dt,\quad \tau \in \mU.
		\label{eq-flqr}
	\end{equation}
	For each $\tau \in \ca U$, let $\beta(\cdot, \tau) \in H^m(\mI)$, the $m$th-order Sobolev space defined by 
	\begin{align*}
		H^m(\mI) = \{h: \mI \to \Re: h \hi {(j)}, j = 0, \ldots, m-1~\mbox{are absolutely continuous, and}~h \hi {(m)} \in L \hi 2 (\mI) \}.
	\end{align*}
	The unknown parameter $\theta_0(\tau) \equiv (\alpha_0(\tau), \beta_0(\cdot, \tau)) \in \ca H \equiv \Re \times H^m(\mI)$.
	We assume $m > 1/2$ to ensure that $H \hi m (\mI)$ is an RKHS.
	
	Given an i.i.d sample, $(X_1, Y_1), \ldots, (X_n, Y_n)$, from the joint distribution of $(X, Y)$, we aim to carry out statistical inference on $\theta_0(\tau)$ in model \eqref{eq-flqr}. The regularized estimator of $(\alpha_0(\tau), \beta_0(\cdot, \tau))$ is defined by 
	\begin{equation}
		(\tilde{\alpha}(\tau), \tilde{\beta}(\cdot, \tau)) = \argmin_{(\alpha, \beta) \in \ca H} \frac{1}{n} \sum_{i = 1}^n \rho_{\tau}\left(Y_i - \alpha - \int_0^1 \beta(t)X(t)dt\right) + \frac{\lambda}{2} J(\beta, \beta) ,
		\label{eq-check}
	\end{equation}
	where $\rho_{\tau}(u) = u[\tau - \bbone(u < 0)]$ is the commonly used loss function in quantile regression (QR), and 
	$J(\beta_1, \beta_2) = \int_0^1 \beta\lo1\hi{(m)}(t) \beta\lo 2 \hi {(m)}(t) dt$ is a roughness penalty for $\beta_1, \beta_2 \in H_m (\mI)$. 
	We use $\lambda/2$ for simplifying future expressions when calculating Fr\'{e}chet derivatives. 
	
	\subsection{Smoothed FLQR}
	As pointed out by \cite{fernandes2019},  the first-order linear Gaussian approximation of the distribution of the standard QR estimator could fail in finite samples. They proposed a convolution-type smoother of the objective function to yield a continuous QR estimator. In particular, the first term of the right-hand side of Equation \eqref{eq-check} can be rewritten as
	\begin{equation}
		\hat{R}(\theta; \tau) := \frac{1}{n}\sum_{i = 1}^n \rho_{\tau}(\epsilon_i(\theta)) = \int \rho_{\tau}(u) d \hat{F}(u; \theta),
		\label{eq-emrisk1}
	\end{equation}
	where $\epsilon_i(\theta) := Y_i - \alpha - \int_0^1 \beta(t)X_i(t)dt$ and $\hat{F}(\cdot; \theta)$ denotes the empirical distribution function of 
	$\epsilon_i(\theta)$'s. Since $\hat{F}(\cdot; \theta)$ is not an absolutely continuous function, the mapping $\theta \to \hat{R}(\theta; \tau)$ is not differentiable. Let $k$ denote a smooth kernel function satisfying $\int k(u) du = 1$ and $h > 0$ be a bandwidth parameter shrinking towards to 0 as the sample size increases. Then define $k_h(v) = \frac{k(v/h)}{h}$. 
	Instead of using the empirical distribution function, \cite{fernandes2019} employed a kernel method to estimate the probability density function of $\eps_i(\theta)$. In our context, the kernel density estimate is given by $\hat{f}(v; \theta) = \frac{1}{n} \sum_{i = 1}^n k_h(v - \eps_i(\theta))$ for any given $\theta \in \ca H$. Hence the corresponding estimated distribution function is $\hat{F}_h(u; \theta) = \int_{-\infty}^u \hat{f}(v; \theta)  dv$. 
	Replacing $\hat{F}(u; \theta)$ with $\hat{F}_h(u; \tau)$ in Equation \eqref{eq-emrisk1} , we obtain the smoothed empirical risk function:
	\begin{equation}
		\hat{R}_h(\theta; \tau) := \int \rho_{\tau}(u) d \hat{F}_h(u; \theta) = \int \rho_{\tau}(u) \hat{f}_h(u; \theta) du. 
		\label{eq-emrisk2}
	\end{equation}
	Then the corresponding regularized FLQR estimator is
	\begin{equation}
		(\hat{\alpha}_{h\lambda}(\tau), \hat{\beta}_{h\lambda}(\cdot, \tau)) \in \argmin_{(\alpha, \beta) \in \ca H} \hat{R}_h(\theta; \tau) + \frac{\lambda}{2}  J(\beta, \beta). 
		\label{eq-smooth}
	\end{equation}

	\subsection{Representer theorem}
	$H^m(\mI)$ is an RKHS when equipped with the (squared) norm
	$$
	\|f\|^2_{\mathcal{W}_2^m} = \sum_{l = 0}^{m - 1} \left\{\int_0^1 f^{(l)}(t) dt\right\}^2 + \int_0^1 \left\{f^{(m)}(t)\right\}^2dt.
	$$
	for any $f \in H^m(\mI)$. Let $H_0(\mI) = \{\beta \in \mH: J(\beta, \beta) = 0\}$ be the null space of $J$, i.e., the collection of all the $\beta$'s such that
	$J(\beta, \beta) = 0$. 
	Then it is a finite-dimensional linear subspace of $H^m(\mI)$. Let $\psi_1, \ldots, \psi_m$ be the basis functions of $H_0(\mI)$. 
	Denote by $H_1(\mI)$ its orthogonal complement in $H^m(\mI)$ such that $H^m(\mI) = H_0(\mI)  \oplus H_1(\mI) $. That is, for any $\beta \in H^m(\mI)$ , there exists a unique decomposition $\beta = \beta_0 + \beta_1$ such that $\beta_0 \in H_0(\mI) $ and $\beta_1 \in H_1(\mI) $. Note that $H_1(\mI) $ is also an RKHS with the inner product of $H^m(\mI)$ restricted to $H_1(\mI)$.

	Let $R$ and $R_1$: $\mI \times \mI \rightarrow \Re$ be the reproducing kernels of $H^m(\mI)$ and $H_1(\mI)$, respectively. Then we have $J(\beta_1, \beta_1) = \|\beta_1\|_R^2 = \|\beta_1\|^2_{\mathcal{W}_2^m} $ for any $\beta_1 \in H_1(\mI)$. According to \cite{cucker2002},
	$(RX_i)(\cdot) := \int_{0}^1 R(\cdot, s)X_i(s)ds \in H^m(\mI)$ for $i = 1, \ldots, n$. As shown in \cite{shin2016rkhs}, the solution to \eqref{eq-smooth} can be expressed as
	\begin{equation}
		\hat{\beta}_{h\lambda}(t, \tau) = \sum_{l = 1}^m d_l\psi_l(t) + \sum_{i = 1}^n c_i\xi_i(t), 
		\label{eq-solution}
	\end{equation}
	where $\xi_i(t) = \int_0^1 R_1(t, s)X_i(s)ds, i = 1, \ldots, n$. Let $\bd = (d_1, \ldots, d_m)^{\trans}, \bc = (c_1, \ldots, c_n)^{\trans}$ and $\bXi = (J(\xi_i, \xi_j))_{ij} \in \Re^{n \times n}$. 
	Therefore,  solving \eqref{eq-smooth} is reduced to minimizing 
	\begin{equation}
		Q_h(\alpha, \bd, \bc; \tau) := \frac{1}{n} \sum_{i = 1}^n \ell_h\left(Y_i - \alpha -  \sum_{l = 1}^m d_l \int_{0}^1 X_i(t)\psi_l(t)dt - \sum_{j = 1}^n c_j \inner{\xi_i}{\xi_j}_{H^m(\mI)} ; \tau \right)      + \frac{\lambda}{2}   \bc \hi {\trans}   \bXi \bc,
		\label{eq-repre}
	\end{equation}
	with respect to $\alpha, \bd$ and $\bc$, where $\ell_h(u; \tau) = \int \rho_{\tau}(v) k_h(v - u)dv$. 
	
	The corresponding estimator of the conditional $\tau$-th quantile of $Y$ is given by $\hat{Q}(\tau, X) = \hat{\alpha}_{h\lambda}(\tau) + \int_0^1 X(t)\hat{\beta}_{h\lambda}(t, \tau)dt$. One possible issue of the above estimation method is that this estimator is not necessarily monotone nondecreasing with respect to $\tau$, a property satisfied by the true conditional quantile of $Y$. 
	Here we adopt the strategy employed in \cite{kato2012} to monotonize the quantile estimator. More specifically, we postulate that $\mU = [\underline{\tau}, \overline{\tau}]$ with $0 < \utau < \otau < 1$. Given $X = x$, we first define a distribution function $\hat{F}_{\mU}(y | x) = \frac{1}{\otau - \utau} \int_{\mU} \bbone \{\hat{Q}_{Y | X}(\tau | x) \leq y\} d\tau$
	with the support $[\min_{\tau \in \mU}\hat{Q}_{Y | X}(\tau | x) , \max_{\tau \in \mU} \hat{Q}_{Y | X}(\tau | x) ]$. Then a modified estimator of ${Q}_{Y | X}(\tau | x)$ is defined by $\tilde{Q}_{Y | X}(\tau | x) = \hat{F}^{-1}_{\mU} ((\tau - \utau)/(\otau - \utau) | x)$. Obviously $\tilde{Q}_{Y | X}(\tau | x) $ is a nondecreasing function of $\tau$. Then we employ the ``pool adjacent violators" algorithm to obtain the isotonized estimator denoted by $\hat{Q}_{Y|X}^I(\tau | x)$. The final estimator of ${Q}_{Y|X}(\tau | x)$ is given by a convex combination of $\tilde{Q}_{Y | X}(\tau | x) $ and $\hat{Q}_{Y|X}^I(\tau | x)$. As shown in \cite{kato2012}, this final estimator of the $\tau$-th quantile, denoted by $\hat{Q}^{\dagger}_{Y | X}(\tau | x)$ is always superior to the initial estimator in the following sense: for any $q \geq 1$, 
	$$
	\left[\int_{\mU} |\hat{Q}^{\dagger}_{Y | X}(\tau | x) - Q_{Y | X}(\tau | x) |^q d\tau \right]^{1/q} \leq \left[\int_{\mU} |\hat{Q}_{Y | X}(\tau | x) - Q_{Y | X}(\tau | x) |^q d\tau \right]^{1/q}.
	$$

	
	\subsection{Gradient descent algorithm}
	To minimize \eqref{eq-repre}, we consider a gradient descent (GD) algorithm. To facilitate calculations of the gradient of the objective function $Q_h$,  we rewrite $\hat{R}_h(\theta; \tau)$ as
	$$
	\hat{R}_h(\theta; \tau) = (1 - \tau)\int_{-\infty}^0 \hat{F}_h(u; \theta) du  + \tau \int_0^{\infty} [1 - \hat{F}_h(u; \theta)] du
	$$
	by the Tonelli's theorem. 
	Then the derivatives of $Q_h(\alpha, \bd, \bc; \tau)$ are immediately available:
	\begin{align*}
		\frac{\partial Q_h(\alpha, \bd, \bc; \tau)}{\partial \alpha} & = \frac{1}{n} \sum_{i = 1}^n\left[\bK\left(-\frac{\ti{\eps}_i(\theta)}{h}\right) - \tau\right],  \\
		\frac{\partial Q_h(\alpha, \bd, \bc; \tau)}{\partial d_l} & = \frac{1}{n} \sum_{i = 1}^n\left[\bK\left(-\frac{\ti{\eps}_i(\theta)}{h}\right) - \tau\right] \left(\int_0^1 X_i(t)\psi_l(t)dt \right), \\
		\frac{\partial Q_h(\alpha, \bd, \bc; \tau)}{\partial \bc} & = \frac{1}{n} \sum_{i = 1}^n\left[\bK\left(-\frac{\ti{\eps}_i(\theta)}{h}\right) - \tau\right]\bXi_i + \lambda \bXi \bc,
	\end{align*}
	where $\ti{\eps}_i(\theta)$ is $\eps_i(\theta)$ with $\beta(t)$ replaced by $\hat{\beta}_h(t, \tau)$ in \eqref{eq-solution},
	$\bXi_i$ is the $i$th column of $\bXi$
	and $\bK(u) = \int_{-\infty}^ u k(v)dv$.  Given these derivatives and an initial value, $(\hat{\alpha}^{0}, \hat{\bd}^{0}, \hat{\bc}^0)$, the GD update at the $(r + 1)$th iteration is given by
	\begin{align}
		\nonumber
		\hat{\alpha}^{r + 1} & = \hat{\alpha}^{r}  - \gamma_r \frac{\partial Q_h(\hat{\alpha}^{r}, \hat{\bd}^{r}, \hat{\bc}^{r}; \tau)}{\partial \alpha}, \\
		\label{eq-GD}
		\hat{d}_l^{r + 1} & = \hat{d}_l^{r}  - \gamma_r \frac{\partial Q_h(\hat{\alpha}^{r}, \hat{\bd}^{r}, \hat{\bc}^{r}; \tau)}{\partial d_l}, ~~l = 1, \ldots, m,  \\
		\nonumber
		\hat{\bc}^{r + 1} & = \hat{\bc}^{r}  - \gamma_r \frac{\partial Q_h(\hat{\alpha}^{r}, \hat{\bd}^{r}, \hat{\bc}^{r}; \tau)}{\partial \bc}, 
	\end{align}
	where $\gamma_r > 0$ is the learning rate. 
	
	Selecting an appropriate learning rate is critical when implementing the GD algorithm. A line search method is commonly used to determine the learning rate. In particular, a fixed learning rate chosen from a tuning strategy or a decreasing sequence of learning rates is taken in practice. However, neither of these approaches is efficient in this context. Instead, we adopt the strategy proposed by \cite{barzilai1988} to select the learning rate at each iteration. For $r = 1, 2, \ldots$, the learning rate $\gamma_r$ is chosen to minimize either $||a\bdel^r -  \bg^r||_2^2$
	or $||\bdel^r - a\bg^r||_2^2$ w.r.t. $a$, where
	\begin{align*}
		\bdel^r & = (\hat{\alpha}^{r}, (\hat{\bd}^{r})^{\trans}, (\hat{\bc}^{r})^{\trans})^{\trans} - (\hat{\alpha}^{r - 1}, (\hat{\bd}^{r - 1})^{\trans}, (\hat{\bc}^{r - 1})^{\trans})^{\trans}  \quad \mbox{and} \\
		\bg^r & = \nabla Q_h(\hat{\alpha}^{r}, \hat{\bd}^{r}, \hat{\bc}^{r}; \tau) - \nabla Q_h(\hat{\alpha}^{r - 1}, \hat{\bd}^{r - 1}, \hat{\bc}^{r - 1}; \tau).
	\end{align*}
	As a result, ordinary least squares lead to two possible choices of the learning rate:
	\begin{equation}
		\gamma_{1, r} = \frac{\inner{\bdel^r}{\bdel^r}}{\inner{\bdel^r}{\bg^r}} \quad\mbox{and}\quad \gamma_{2, r} = \frac{\inner{\bdel^r}{\bg^r}}{\inner{\bg^r}{\bg^r}}.
		\label{eq-BBstep}
	\end{equation}
	Both of them are referred to as the BB rate in the article. Then the GD update proceeds as
	in \eqref{eq-GD} by replacing $\gamma_r$ with either one in \eqref{eq-BBstep}. 
	Note that we need point estimates of all parameters at the first iteration to calculate $\gamma_{1, 1}$ and/or $\gamma_{1, 2}$, but we only have an initial value ($r = 0$) at hand. Thus we suggest taking $\gamma_0$ to be 1 in \eqref{eq-GD} to obtain the point estimates at the first iteration. 
	
	As pointed out in \cite{he2020}, the BB rate could be negative at some iterations and may vibrate drastically sometimes especially when $\tau$ approaches to 0 or 1. To address these issues, we follow their suggestion to choose $\gamma_r = \min\{\gamma_{1,r}, \gamma_{2, r}, 100\}$ if $\gamma_{1, r} > 0$ and 1 otherwise. 
	Iterations of our proposed GD-BB algorithm will not be stopped until $||\nabla Q_h(\hat{\alpha}, \hat{\bd}, \hat{\bc}; \tau)||_2 \leq tol$, where $tol > 0$ is a predetermined threshold called the gradient tolerance in literature.

	\section{Theoretical Properties}
	\label{sec:Bahadur}
	Let $R(\theta; \tau)$ and $R_h(\theta; \tau)$ denote the expected loss function in \eqref{eq-emrisk1} and \eqref{eq-emrisk2}, respectively. Note that $\hat{\theta}_{h\lambda}(\tau)$ actually estimates $\theta_{0h}(\tau) = \argmin_{(\alpha, \beta) \in \ca H} R_h(\theta; \tau)$. 
	To investigate the consistency and weak convergence of $\hat{\theta}_{h\lambda}(\tau)$, we introduce the following objective function 
	\begin{equation} \label{eq-PPL}
		\ell_{h\lambda} := R_h(\theta; \tau) + \frac{\lambda}{2} J(\beta, \beta) 
	\end{equation}
	to bridge the gap between $\hat{\theta}_{h\lambda}(\tau)$ and $\theta_0(\tau)$. Our main idea is to first establish the approximation error of the minimizer of \eqref{eq-PPL}, denoted by $\theta_{0h,\lambda}(\tau)$, relative to $\theta_0(\tau)$, and then quantify the estimation error of $\hat{\theta}_{h\lambda}(\tau)$ relative to $\theta_{0h,\lambda}(\tau)$. We develop a key technical tool called the {\it functional Bahadur representation} of the estimator $\hat{\theta}_{h\lambda}(\tau)$  to derive its consistency and weak convergence. 
	
	%
	%
	%

	\subsection{Approximation error analysis}
	
	Let $f (\cdot | x)$ be the conditional probability density function of $Y$ given $X = x$. For the existence of the regular conditional distribution of $Y$ given $X$, please refer to \cite{kato2012} for a more detailed discussion. To establish the desirable theoretical properties of the proposed estimator, we need the following regularity conditions. 
	\begin{assump} \label{ass-f}
		The conditional density and quantile of $Y$ given $X = x$ satisfy
		
		\begin{itemize}
			
			
			\item [(a)] The conditional density $f(y | x)$ is continuous and strictly positive over $\Re \times \mbox{supp}(X)$. 
			
			\item [(b)] There exists an integer $s \geq 1$ such that the $s$-th order derivative $f \hi {(s)}(\cdot | x)$ is uniformly continuous in the sense that 
			$$
			\lim\limits_{\epsilon \to 0} \sup_{(y, x) \in \Re \times \mbox{supp}(X)} \sup_{t: |t| \leq \epsilon}|f^{(s)}(y + t | x)
			- f^{(s)}(y | x)| = 0,
			$$
			and that for $j = 0, 1\ldots,s$, $\sup_{(y, x) \in \Re \times \mbox{supp}(X)}  |f^{(j)}(y | x)| < \infty$ and $\lim_{y \to  \pm\infty}f^{(j)}(y | x) = 0$. 
			
			\item [(c)] $X$ is weighted-centered in the sense that $\E\{B_{\tau}(X)X(t)\} = 0$ for every $t \in \mI$, where \\
			$B_{\tau}(X) = f\left(\alpha_{0}(\tau) + \int_0^1 \beta_{0}(t, \tau)X(t) dt \big| X \right)$. Furthermore, we assume that there exists a positive constant $C_1$ such that $C_1^{-1} \leq B_{\tau}(X) \leq C_1\quad\mbox{a.s.}$
			
		\end{itemize}
	\end{assump}
	\begin{remark} \label{re-B(X)}
		Assumptions \ref{ass-f} (a) and (b) are fairly common in the literature of quantile regression; see \cite{belloni2019}, \cite{chao2017quantile} and  \cite{fernandes2019} for instance. 
		Concerning the definition of $B_{\tau}(X)$, it is actually $f_{\epsilon}(F^{-1}_{\eps}(\tau))$ in the context of a standard functional linear regression $Y = \alpha + \int \beta(t)X(t)dt + \epsilon$, where $\eps$ is independent of $X$. 
		
	\end{remark}	
	
	\begin{assump} \label{ass-K}
		The kernel function $k$ and bandwidth $h$ satisfy
		\begin{itemize}
			\item [(a)]
			The kernel function $k: \Re \to \Re$ is even, integrable, twice differentiable with bounded first and second derivatives, and satisfies that $\int k(u)du = 1$ and $0 < \int_0^{\infty} \bK(u)[1 - \bK(u)]du < \infty$. Additionally, $\int |u^{s + 1}k(u)| du < \infty$ for $s$ in (b) of Assumption \ref{ass-f}, and $k$ is orthogonal to all non-constant monomials of degree up to $s$, i.e., $\int u^j k(u)du = 0$ for $j = 1, \ldots, s$, and $\int u^{s + 1}k(u) du \neq 0$. 
			
			\item [(b)] $h = h_n = o(1)$ as $n \rightarrow \infty$.

		\end{itemize}
	\end{assump}

	Let $C(s, t; \tau) = \E[B_{\tau}(X)X(t)X(s)]$ be the weighted covariance function. We introduce the following inner product in $H^m(\mI)$:
	\begin{equation} \label{eq-inner1}
		\inner{\beta}{\tilde{\beta}}_{1} = V(\beta, \tilde{\beta}) + \lambda J(\beta, \tilde{\beta}),
	\end{equation}
	where $V(\beta, \tilde{\beta}) \equiv \int_0^1 \int_0^1 C(s, t; \tau) \beta(t)\tilde{\beta}(s) dsdt$. Here we suppress the dependence of $V$ on $\tau$ for ease of notation. 
	Define an integral operator $C$ from $L^2(\mI)$ to $L^2(\mI)$: $\beta \to (C\beta)(t) = \int_0^1 C(s, t; \tau)\beta(s)ds$. 
	The following regularity condition ensures that $V$ is positive definite.

	\begin{assump}
		
		$C(s, t; \tau)$ is continuous on $\mI \times \mI$. Furthermore, for any $\beta \in L^2(\mI)$, $C\beta = 0$ if and only if $\beta = 0$. 
		\label{ass-positive}
	\end{assump}
	
	As illustrated in \cite{shang2015}, $H^m(\mI)$ is an RKHS under $\inner{\cdot}{\cdot}_1$. Denote the reproducing kernel function by $K(s, t)$. Then we define a proper inner product in $\ca H$ such that it is a well-defined Hilbert space. Following the idea of \cite{shang2015}, for any $\theta = (\alpha, \beta), \ti{\theta} = (\ti{\alpha}, \ti{\beta}) \in \ca H$, define 
	\begin{equation} \label{eq-inner2}
		\inner{\theta}{\ti{\theta}} = \E\left\{B_{\tau}(X) \left(\alpha + \int_0^1 X(t)\beta(t)dt\right) \left(\ti\alpha + \int_0^1 X(t)\ti\beta(t)dt\right) \right\}+ \lambda J(\beta, \ti{\beta}).
	\end{equation}
	Obviously $\inner{\theta}{\ti{\theta}}  =  \inner{\beta}{\tilde{\beta}}_1 + \E[B_{\tau}(X)]\alpha\ti{\alpha}$. Given the fact that $V$ is positive definite under Assumption \ref{ass-positive}, $\inner{\cdot}{\cdot}$ defined in \eqref{eq-inner2} is a well-defined inner product. This indicates that $\ca H$ is indeed a Hilbert space equipped with the norm $\inner{\cdot}{\cdot}$ in \eqref{eq-inner2}.
	


	Next we assume that there exists a sequence of basis functions in $H^m(\mI)$ that can simultaneously diagonalize both $V$ and $J$ in Equation \eqref{eq-inner1}. 
	
	\begin{assump}
		\label{ass-diag}
		There exists a sequence of basis functions $\{\phi_{\nu}\}_{\nu \geq 1} \subsetneq H^m(\mI)$ such that 
		$\|\phi_{\nu}\| \lo {L^2} \leq C_{\phi} \nu^a$ holds uniformly over $\nu$ for some constants $a \geq 0$, $C_{\phi} > 0$, and that
		\begin{equation}
			V(\phi_{\mu}, \phi_{\nu}) = \delta \lo {\mu\nu}, \quad J(\phi_{\mu}, \phi_{\nu}) = \rho\lo{\nu}\delta \lo {\mu\nu}
			\quad\mbox{for any~} \mu, \nu \geq 1. 
			\label{eq-diag}
		\end{equation}
		where $\delta_{\mu\nu} = 1$ if $\mu = \nu$ and 0 otherwise, and $\rho_{\nu}$ is a nondecreasing nonnegative sequence satisfying $\rho \lo {\nu} \asymp \nu^{2l}$ for some constant $l > a + 1/2$. Furthermore, any $\beta \in H^m(\mI)$ admits the Fourier expansion $\beta = \sum_{\nu = 1}^{\infty} V(\beta, \phi_{\nu}) \phi \lo {\nu}$ with convergence in $H^m(\mI)$ under $\inner{\cdot}{\cdot}_1$. 
	\end{assump}
	
	Let $K_t(\cdot) := K(t, \cdot) \in H^m(\mI)$ for any $t \in \mI$. Assumption \ref{ass-diag} implies that there exists a sequence of real numbers $a_{\nu}$ such that $K_t = \sum_{\nu = 1}^{\infty} a_{\nu} \phi \lo {\nu}$. Obviously, 
	$\phi\lo{\nu}(t) = \inner {K_t} {\phi \lo {\nu}}_1 = a_{\nu}(1 + \lambda \rho \lo {\nu})$ for any $\nu \geq 1$. It follows that
	$K_t = \sum_{\nu = 1}^{\infty} \frac{\phi \lo {\nu} (t)}{1 + \lambda \rho \lo {\nu}} \phi \lo {\nu}$. By the Riesz representation theorem, there exists  a linear operator $W_{\lambda}$ from $H^m(\mI)$ to $H^m(\mI)$ such that
	$\inner{W\lo{\lambda} \beta} {\ti{\beta}}_1 = \lambda J(\beta, \ti{\beta})$ for all $\beta, \ti{\beta} \in H^m(\mI)$. By the definition of $J$, $W\lo{\lambda}$ is nonnegative and self-adjoint. For any $\nu \ge 1$, write $W\lo{\lambda} \phi \lo {\nu} = \sum_{\mu = 1}^{\infty} b\lo{\mu} \phi \lo {\mu}$. From Assumption \ref{ass-diag}, we have 
	$\lambda \rho \lo{\nu} \delta \lo {\nu\mu} = \lambda J(\phi \lo {\nu}, \phi \lo {\mu}) = \inner{W\lo{\lambda} \phi \lo {\nu} }{\phi \lo {\mu}}_1 = b_{\mu}(1 + \lambda \rho_{\mu}) $. Therefore, $b_{\mu} = \lambda \rho_{\nu} /(1 + \lambda \rho_{\nu})$ if $\mu = \nu$ and 0 otherwise. As a result, $W\lo{\lambda} \phi \lo {\nu} =  \frac{\lambda \rho_{\nu}}{1 + \lambda \rho_{\nu}}\phi \lo {\nu}$
	
	\begin{proposition} \label{propW}
		Suppose that Assumption \ref{ass-diag} holds, then for any $t \in \mI$, we have 
		$$
		K_t(\cdot) = \sum_{\nu = 1}^{\infty} \frac{\phi \lo {\nu} (t)}{1 + \lambda \rho \lo {\nu}} \phi \lo {\nu}(\cdot),
		$$
		and for any $\nu \ge 1$, 
		$$
		(W\lo{\lambda} \phi \lo {\nu})(\cdot) =  \frac{\lambda \rho_{\nu}}{1 + \lambda \rho_{\nu}}\phi \lo {\nu}(\cdot). 
		$$
	\end{proposition}
	For any $x \in L^2(\mI)$, let $\eta(x) = \sum_{\nu = 1}^{\infty} \frac{x \lo {\nu}}{1 + \lambda \rho \lo {\nu}} \phi \lo {\nu}$, where
	$x \lo {\nu} = \int_0^1 x(t) \phi \lo {\nu}(t) dt$.
	By Proposition 2.4 of \cite{shang2015}, we can conclude that $R_x = ([\E\{B_{\tau}(X)\}]^{-1}, \eta(x)) \in \ca H$, satisfies
	$\inner{R_x}{\theta} = \alpha + \int_0^1 x(t)\beta(t)dt$ for any $\theta = (\alpha, \beta) \in \ca H$. We need to emphasize that $R_x$ depends on $\lambda$ through the definition of $\eta(x)$. 
	
	Let $\ell_{n, h\lambda}(\theta) = \hat{R}_h(\theta; \tau)  + \frac{\lambda}{2} J(\beta, \beta)$. For notational convenience, denote  
	$\Delta \theta = (\Delta \alpha, \Delta \beta)$ and $\Delta \theta_j = (\Delta \alpha_j, \Delta \beta_j)$ for $j = 1, 2, 3$. 
	The first order Fr\'{e}chet derivative operator of $\ell_{n, h\lambda}(\theta) $ w.r.t $\theta$ is given by
	\begin{equation}
		S_{n,h\lambda}(\theta) \Delta \theta  =    \frac{1}{n} \sum_{i = 1}^n\left[\bK\left(\frac{-\eps_i(\theta)}{h}\right) - \tau\right]\inner{R \lo {X_i}}{\Delta \theta}+ \lambda J(\beta, \Delta \beta).
		\label{eq-first}
	\end{equation}
	The second- and third-order Fr\'{e}chet derivatives of $\ell_{n, h\lambda}$ can be shown to be, respectively, 
	$$
	DS_{n, h\lambda}(\theta)\Delta \theta_1 \Delta \theta_2 =  \frac{1}{n} \sum_{i = 1}^nk_h(-\eps_i(\theta))\inner{R \lo {X_i}}{\Delta \theta_1} \inner{R \lo {X_i}}{\Delta \theta_2}+ \lambda J(\Delta\beta_1, \Delta \beta_2)
	$$
	and 
	$$
	D^2S_{n, h\lambda}(\theta)\Delta \theta_1 \Delta \theta_2 \Delta \theta_3 =  \frac{1}{n} \sum_{i = 1}^nh^{-2}k^{(1)}(-\eps_i(\theta)/h)\inner{R \lo {X_i}}{\Delta \theta_1} \inner{R \lo {X_i}}{\Delta \theta_2} \inner{R \lo {X_i}}{\Delta \theta_3} 
	$$
	
	
	\begin{assump} 	\label{ass-X}
		The functional covariate $X$ satisfies 	that there exists a constant $s \in (0, 1)$ such that

		\begin{equation}
			\E\{\exp(s\|X\|_{L^2})\} < \infty,
			\label{eq-mgf}
		\end{equation}
		where $\|X\|_{L^2}^2 = \int_0^1 X^2(t)dt$. We further assume that there exists a constant $M_0 > 0$ such that for any $\beta \in H^m(\mI)$, 
		\begin{equation}
			\E\left\{\left|\int_0^1 X(t)\beta(t)dt\right|^4\right\} \le M_0 \left[   \E\left\{\left|\int_0^1 X(t)\beta(t)dt\right|^2\right\}              \right]^2
			\label{eq-moment}
		\end{equation}
		
	\end{assump}
	
	\begin{remark} \label{re-subexp}
		It is easy to check that \eqref{eq-mgf} holds if the support of $X$ is a compact subset of $L^2(\mI)$ since $\|X\|_{L^2} \le c_0$ a.s. for some constant $c_0 > 0$. 
		When $X$ is a Gaussian process with a square-integrable mean function, Equation \eqref{eq-mgf} holds for any $s \in (0, 1/4)$ by Proposition 3.2 of \cite{shang2015},
		and Equation \eqref{eq-moment} is fulfilled for $M_0 = 3$; see \cite{yuan2010}. 
	\end{remark}

	For $j = 1, 2, 3$, let $R^{(j)}
	(\theta; \tau)$ and $R_h^{(j)}(\theta; \tau)$ be the $j$th order Fr\'{e}chet derivative operators of $R(\theta; \tau)$ and $R_h(\theta; \tau)$, respectively.
	Denote by $\ca S$ the set $\ca H \times [\underline{\tau}, \bar{\tau}]$ to which $(\theta, \tau)$ belongs.  An auxiliary norm is introduced for technical purpose: $\|\theta\|_2 = |\alpha| + \|\beta\|\lo{L^2}$ for any $\theta =(\alpha, \beta) \in \ca H$. Before we present the theorem concerning the approximation error of $\theta_{0h, \lambda}(\tau)$, we first establish the following properties of the expected surrogate loss function $R_h(\theta; \tau)$. 
	\begin{lemma}
		\label{lemma1}
		Suppose Assumptions \ref{ass-f}, \ref{ass-K} and \ref{ass-X} hold. Then
		\begin{itemize}
			\item [(i)] $\sup_{(\theta, \tau) \in \ca S} \left| \frac{R_h(\theta; \tau) - R(\theta; \tau)}{h^{s+1}} \right| = O(1)$;
			\item [(ii)] $\sup_{(\theta, \tau) \in {S}} \left| \frac{R_h^{(1)}(\theta; \tau)\theta_1 - R^{(1)}(\theta; \tau)\theta_1}{h^{s+1}} \right| = O(1)\|\theta_1\|$ for any $\theta_1 \in \ca H$. 
			\item [(iii)] $\sup_{(\theta, \tau) \in {S}} \left| \frac{R_h^{(2)}(\theta; \tau)\theta_1\theta_2 - R^{(2)}(\theta; \tau)\theta_1\theta_2}{h^{s}} \right| = o(1)\|{\theta_1}\|\|{\theta_2}\|$ for any $\theta_1, \theta_2 \in \ca H$. 
			\item [(iv)] $\sup_{(\theta^*, \theta, \tau) \in {\ca H \times S}} \left| \frac{R_h^{(2)}(\theta + \theta^*; \tau)\theta_1\theta_2 - R_h^{(2)}(\theta; \tau)\theta_1\theta_2}{\|\theta^*\|} \right| = O(1)\|\theta_1\|\|\theta_2\|$ for any $\theta_1, \theta_2 \in \ca H$. 
		\end{itemize}
	\end{lemma}

	The following theorem establishes the approximation error of the minimizer of $\ell_{h\lambda}$ relative to $\theta_0(\tau)$.
	
	\begin{theorem} \label{th-thlam}
		Suppose that Assumptions \ref{ass-f} to \ref{ass-X} hold, and there exists some constant $C_2 > 0$ such that $\sup_{\tau \in [\underline{\tau}, \bar{\tau}] } J(\beta_0(\cdot, \tau), \beta_0(\cdot, \tau)) \leq C_2$ 
		Then there exists a unique minimizer of $\ell_{h \lambda}$, denoted by $\theta_{0h, \lambda}(\tau)$, for every $\tau \in [\underline{\tau}, \bar{\tau}]$. Furthermore, $\|\theta_{0h, \lambda}(\tau) - \theta_0(\tau)\| = O(\lambda^{1/2})$ uniformly over $\tau \in [\underline{\tau}, \bar{\tau}]$, provided that $h = O(\lambda^{1/(2s+2)})$. 
		
	\end{theorem}

	\subsection{Bahadur representation for smoothed FLQR}
	
	Next we study the estimation error of $\hat{\theta}_{h\lambda}(\tau)$. Denote $\omega = \lambda^{1/(2l)}$, where $l$ is specified in Assumption \ref{ass-diag}. Before laying out the main theorem, we first prove that $\hat{\theta}_{h\lambda}(\tau)$ is a consistent estimator of $\theta_0(\tau)$ and the convergence rate is established.

	\begin{theorem} [Convergence rate]
		Suppose that Assumptions \ref{ass-f} to \ref{ass-X} hold, and that the following rate conditions on $\omega$ (or equivalently, $\lambda$) and $h$ are satisfied 
		\begin{align}
			\nonumber
			h = O(w^{l/(s+1)}),\quad\quad	\omega  = o(1), \quad\quad n^{-1/2}\omega^{-1}h^{-2} & = o(1), \quad\quad	\\
			\label{eq-lamorder}
			n^{-1/2}\omega^{-(a+1) - ((2l-2a-1)/(4m))}  (h^{-1}\log n) (\log \log n)^{1/2} &  = o(1). 
		\end{align}
		Then $\hat{\theta}_{h\lambda}(\tau) = (\hat{\alpha}_h(\tau), \hat{\beta}_h(\cdot,\tau))$ is the unique solution to \eqref{eq-smooth} and for every $\tau \in [\underline{\tau}, \bar{\tau}]$, $||\hat{\theta}_{h\lambda}(\tau)  - \theta_{0}(\tau)|| = O_P(r_n)$, where $r_n = (n\omega)^{-1/2} + \omega^l$. 
		\label{th-rate}
	\end{theorem}

	Now we are ready to present our main theorem: the Bahadur representation for the FLQR model under the RKHS framework. Similar results have been established in classic semi-/non-parametric regression; see \cite{shang2010ejs}, \cite{shang2013aos}, and \cite{shang2015aos}.
	Our result is the first one in functional quantile regression. 
	
	\begin{theorem} [Functional Bahadur representation] 		\label{th-Bahadur}
		Suppose that the assumptions of Theorem \ref{th-rate} are met.  
		Then  as $n \rightarrow \infty$, $\|\hat{\theta}_{h,\lambda}(\tau) - \theta_0(\tau) + S_{n, h\lambda}(\theta_0(\tau)) \| = O_P(a_n)$ for every $\tau \in [\underline{\tau}, \bar{\tau}]$, where
		$$
		a_n = n^{-1/2} \omega^{-\frac{4ma+6m-1}{4m}} r_n( h^{-1} \log n) (\log \log n)^{1/2} + 2h^{-1}r_n^2
		$$
	\end{theorem}

	The Bahadur representation implies the following pointwise limiting distribution of the estimated slope function defined by Equation \eqref{eq-smooth}. This result will be used to construct pointwise confidence intervals for $\beta_0(\cdot, \tau)$.

	\begin{corollary} \label{cor-pointwise}
		Suppose that the assumptions of Theorem \ref{th-Bahadur} are satisfied, and $\sup_{\nu \geq 1} \|\phi_{\nu}\| \leq C_{\phi}\nu^a$ where the constants $a$ and $C_{\phi}$ are defined in Assumption \ref{ass-diag}. 
		Furthermore, we assume $\sum_{\nu} \frac {\phi_{\nu}^2(t)} {(1 + \lambda\rho_{\nu})^2} \asymp \omega^{-(2a + 1)}$, $n^{1/2}a_n = o(1)$, $n^{1/2}h^{s+1} = o(1)$ and $n\omega^{2a+1}(\log(\omega^{-1}))^{-4} \rightarrow \infty$, as $n \rightarrow \infty$. Let $s_n^2 = \Var\{\sum_{i = 1}^n [\bK({-\eps_i(\theta_0(\tau))}/{h}) - \tau] \eta(X_i)(t)\}$.
		Then we have for any fixed $t \in \mI$ and any fixed $\tau \in \mU$, 
		$$
		\frac{n}{s_n}[\hat{\beta}_{h,\lambda}(t, \tau) - \beta_0(t, \tau) + (W_{\lambda}\beta_0(\cdot, \tau))(t)]  \convd N(0, 1).
		$$
		Additionally, if $ns_n^{-1} (W_{\lambda}\beta_0(\cdot, \tau))(t) = o(1)$, then 
		$$
		\frac{n}{s_n}[\hat{\beta}_{h,\lambda}(t, \tau) - \beta_0(t, \tau)]  \convd N(0, 1).
		$$
		
	\end{corollary}

	\begin{remark} \label{re-pointCI}
		As explained in Remark \ref{re-B(X)}, if the underlying true model is given by $Y_i = \alpha + \int \beta(t)X_i(t)dt + \epsilon$, where $\eps$ is independent of $X$, $B_{\tau}(X) = f_{\eps}(F^{-1}_{\eps}(\tau))$. In this scenario, 
		based on the proof Corollary \ref{cor-pointwise}, as $n \rightarrow \infty$, $s_n^2$ can be taken as $\frac{n\tau(1 - \tau)}{f_{\eps}(F^{-1}_{\eps}(\tau))} \sum_{\nu} \frac {\phi_{\nu}^2(t)} {(1 + \lambda\rho_{\nu})^2}$. Then given any fixed $t_0 \in \mI$ and $\tau \in \mU$, a $(1 - \xi)$ pointwise confidence interval for $\beta_0(t_0, \tau)$  is immediately available:
		$$
		P\left( \beta_0(t_0, \tau) \in \left[\hat{\beta}_{h,\lambda}(t_0, \tau) +  \pm z_{\xi/2} \sqrt{\frac{{\tau(1-\tau)}}{{nf_{\eps}(F^{-1}_{\eps}(\tau))}} \cdot \sum_{\nu}\frac {\phi_{\nu}^2(t_0)} {(1 + \lambda\rho_{\nu})^2}}~\right]  \right) \rightarrow 1 - \xi
		$$
		as $n \rightarrow \infty$, where $z_{\xi/2}$ is the $(1 - \xi/2)$-quantile of $N(0, 1)$, if $ns_n^{-1}[W_{\lambda}\beta_0(\cdot, \tau)](t_0) = o(1)$.

	\end{remark}
	
	Besides its limiting distribution at a fixed point $t \in \mI$, we can even establish weak convergence of $\hat{\beta}_{h,\lambda}(\cdot, \tau)$ in the Hilbert space $H^m(\mI)$ equipped with the inner product
	$V(\cdot, \cdot)$ for any given $\tau \in \mU$. Define $\kappa_{\nu}^2 := (\tau - \tau^2)\E[\int_0^1 X_i(t)\phi_{\nu}(t)dt]^2$. 
	
	\begin{theorem} \label{th-weak}
		Suppose that the assumptions of Theorem \ref{th-Bahadur} are satisfied, and $\sup_{\nu \geq 1} \|\phi_{\nu}\| \leq C_{\phi}\nu^a$ where the constants $a$ and $C_{\phi}$ are defined in Assumption \ref{ass-diag}. 
		Furthermore, $\sum_{\nu} \frac {\phi_{\nu}^2(t)} {(1 + \lambda\rho_{\nu})^2} \asymp \omega^{-(2a + 1)}$ holds uniformly for $t \in \mI$, $n^{1/2}a_n = o(1)$, $n^{1/2}h^{s+1} = o(1)$, $n^{1/2}\lambda = o(1)$ and $n\omega^{2a+1}(\log(\omega^{-1}))^{-4} \rightarrow \infty$, as $n \rightarrow \infty$. 
		In addition, the true slope function $\beta_0(\cdot, \tau)$ can be written as $\beta_0 = \sum_{\nu} b_{\nu}\phi_{\nu}$ satisfying $\sum\lo{\nu} b_{\nu}^2 \rho_{\nu}^2 < \infty$,
		and $\Gamma(s, t) = \lim_{n\to\infty}\omega^{2a + 1} \sum_{\nu} \frac{\kappa_{\nu}^2 \phi_{\nu}(s) \phi_{\nu}(t)}{(1 + \lambda\rho_{\nu})^2}$ finitely exists for any $s,t\in[0,1]$.
		Then $n^{1/2}\omega^{a + 1/2}[\hat{\beta}_{h,\lambda}(t, \tau) - \beta_0(t, \tau)]$ converges weakly to a mean zero Gaussian process $\ca G (t)$ in the Hilbert space $H^m(\mI)$ equipped with the inner product $V(\cdot,\cdot)$, whose covariance function is $\Gamma(s,t)$.
	\end{theorem}

	\begin{remark} \label{re-SCB}
		We leverage Theorem \ref{th-weak} to construct an $(1 - \alpha)$ simultaneous confidence band (SCB) for $\beta(\cdot, \tau)$.
		In particular, we need to find a critical value $q_{\alpha}$ that approximately satisfies 
		$
		P\left(\sup_{t \in \mI} | \ca G (t)| \leq q_{\alpha} \right) = 1 - \alpha.
		$
		We employ the \textbf{Chebfun} package in Matlab to find the eigenvalues $\rho_{\nu}$'s and the corresponding eigenfunctions $\phi_{\nu}$'s by solving the integro-differential
		equations (2.7) of \cite{shang2015}. One can refer to \cite{driscoll2008} for more details on the implementation of this open source package. 
		After estimating the covariance function $\Gamma(s, t)$ of $\ca G(t)$. we generate independent sample paths of $\ca G(t)$ to estimate the critical value $q_{\alpha}$. 

	\end{remark}

	\subsection{Confidence interval for conditional quantile} \label{sec-CIquant}
	In this section, we discuss confidence intervals for conditional quantile. In particular, given any nonrandom $x_0 \in L^2(\mI)$, we aim to construct a confidence interval for the true $\tau$-th conditional quantile of $Y$ given $X = x_0$. Obviously $Q_{Y | X}(\tau | x_0) = \alpha_0(\tau) + \int_0^1 x_0(t)\beta_0(t, \tau)dt$. A direct plug-in point estimate of $Q_{Y | X}(\tau | x_0)$ is $\hat{Q}_{Y | X}(\tau | x_0) = \hat{\alpha}_{h,\lambda}(\tau) + \int_0^1 x_0(t) \hat{\beta}_{h,\lambda}(t, \tau)dt$. Let $\sigma_n^2(x_0) = [\E\{B_{\tau}(X)\}]^{-1} + \sum_{\nu} \frac{|x_{\nu}^0|^2}{(1 + \lambda\rho_{\nu})^2}$, where $x_{\nu}^0 = \int_0^1 x_0(t) \phi_{\nu}(t) dt$, $\nu \geq 1$. The following theorem provides a valid $(1 - \xi)$ confidence interval for the true conditional quantile.  
	
	\begin{theorem} \label{th-quantCI}
		Suppose that the assumptions of Theorem \ref{th-Bahadur} are satisfied. 
		Furthermore, $n^{1/2}a_n = o(1)$, $\|R_{x_0}\| \asymp \sigma_n(x_0)$, $n^{1/2}h^{s+1} = o(1)$ and $n\omega^{2a+1}(\log(\omega^{-1}))^{-4} \rightarrow \infty$, as $n \rightarrow \infty$. Let $s_{2n}^2 = \Var\left(\sum_{i = 1}^n [\bK\{{-\eps_i(\theta_0(\tau))}/{h}\} - \tau] \inner{R_{x_0}}{R_{X_i}}\right)$.
		Then
		
		$$
		\frac{n}{s_{2n}} \left[\hat{Q}_{Y|X}(\tau | x_0) - Q_{Y | X}(\tau | x_0) + \int_0^1 x_0(t) \{W_{\lambda}(\beta_0(\cdot,\tau))\}(t)dt \right]  \convd N(0, 1).
		$$
		Additionally, if the true slope function $\beta_0(\cdot, \tau)$ can be written as $\beta_0(\cdot, \tau) = \sum_{\nu} b_{\nu}\phi_{\nu}$ satisfying $\sum\lo{\nu} b_{\nu}^2 \rho_{\nu}^2 < \infty$ and $n\omega^{4l} = o(1)$, then $\frac{n}{s_{2n}} \int_0^1 x_0(t) (W_{\lambda}(\beta_0))(t)dt = o(1)$. It follows that  
		$$
		\frac{n}{s_{2n}}[\hat{Q}_{Y|X}(\tau | x_0) - Q_{Y | X}(\tau | x_0)]  \convd N(0, 1).
		$$
		Hence the $(1 - \xi)$ confidence interval for $Q_{Y | X}(\tau | x_0)$ is 
		$$
		[\hat{Q}_{Y|X}(\tau | x_0) \pm n^{-1}s_{2n} z_{\xi/2}]
		$$
		
	\end{theorem}

	\section{Simulation study}
	\label{sec:simulation}
	
	In this section, we carry out simulation studies to evaluate the
	finite-sample performance of the proposed method in estimating the coefficient function in FLQR. We further investigate the estimated confidence intervals and the SCB proposed in Section \ref{sec:Bahadur}. 
	
	\subsection{Estimation of $\beta(\cdot, \tau)$}  
	
\begin{figure}[H]
		\centering
		{\includegraphics[width=8cm]{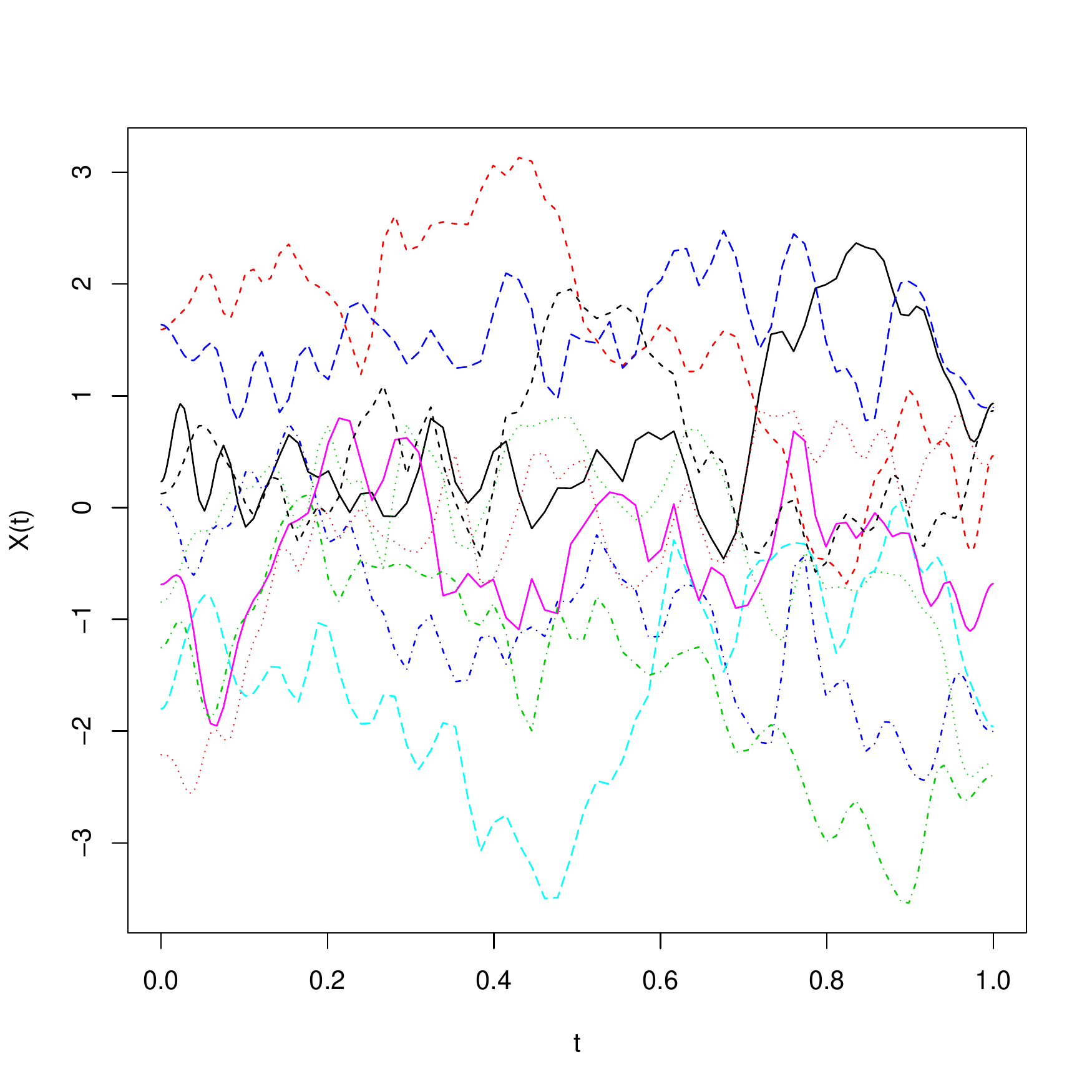}}
		\caption{Trajectories of the functional covariate from 10 randomly selected subjects in one simulation run.}
		\label{fig:design}
	\end{figure}

	We first compare the performance of the estimated coefficient function through the FPCA approach proposed by \cite{kato2012} and our proposed approach under the RKHS framework. In particular, the functional covariate is generated as $X_i(t) = \sum_{k = 1}^{50} \zeta_k \xi_{ik} \psi_k$ for $i = 1, \ldots, n$, where $\zeta_k = 4\cdot(-1)^{k+1} \cdot k^{-2}$, and $\psi_1 = 1$ and $\psi_k = \sqrt{2}\cos((k-1)\pi t)$ for $k \geq 2$ are the eigenfunctions of $X$. The random coefficients $\xi_{ik}$'s are i.i.d generated from uniform distribution on $[-\sqrt{3}, \sqrt{3}]$.
	Figure \ref{fig:design} displays the trajectories of 10 randomly selected functional covariates. 
	The response is generated as $y_i = \alpha + \int_0^1 X_i(t)\beta(t) dt + \sigma\epsilon_i$, where the intercept $\alpha = 0.1$ and the coefficient function $\beta(t) = e^{-t}$. We consider two different designs, normal distribution of mean 0 and student $t$ distribution with 3 degrees of freedom, for the random noise $\epsilon$ to accommodate both light and heavy tails. Different values of $\sigma$'s are chosen to deliver two signal-to-noise ratios (SNR): 10 and 5. 200 independent Monte Carlo simulations are run for each design, where $n = 200$ curves are generated.
	
	\begin{figure}[H]
		\centering
		{\includegraphics[width=12cm]{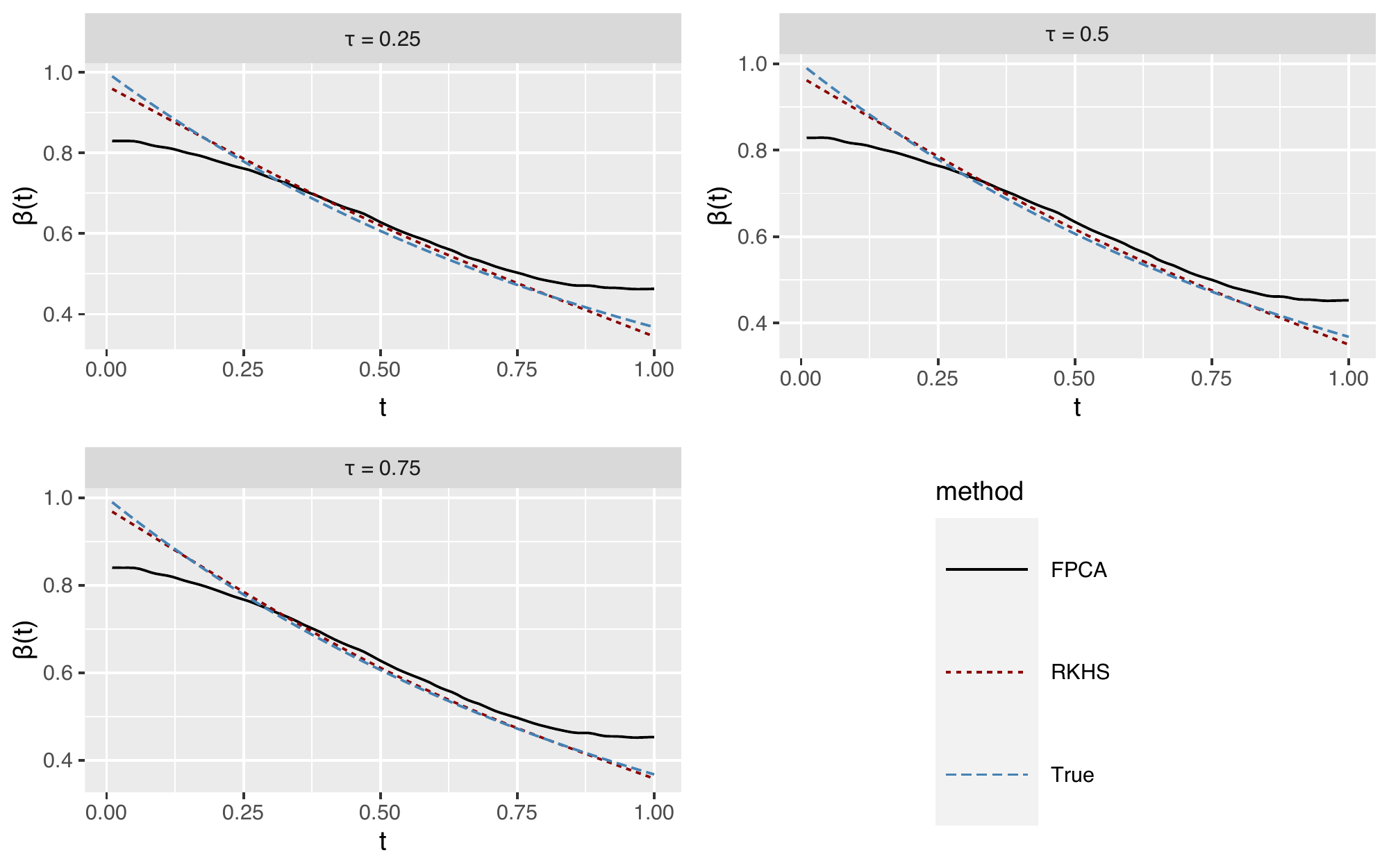}}
		\caption{True slope function $\beta(t)$ and estimates from FPAC and RKHS based methods. In each panel, the solid black line represents the estimate from the FPCA-based method, while the dashed blue line and the dotted red line represent the true coefficient function and the estimate from the proposed method, respectively.}
		\label{fig:estimate}
	\end{figure}
	
	To implement our proposed method, we take $m = 2$ for the choice of $H^m(\mI)$ and adopt the \cite{silverman1986}'s rule-of-thumb (ROT) bandwidth to choose $h$ in kernel smoothing for saving computational cost. More specifically, we first fit a standard quantile regression of $Y$ with respect to $z_{ij} = \int X_i(t)g_j(t) dt$, where $g_1 = 1, g_2(t) = t - 0.5$ and $g_j(t) = \xi_{j-2}(t)$ for $j = 3, \ldots, 2+n$.  Then the ROT bandwidth is given by $h_{\text{ROT}} = 1.06\hat{s}n^{-1/5}$, where $\hat{s}$ is the minimum of the
	estimated standard error and the interquartile range (divided
	by 1.39) of the residuals from fitting the standard quantile regression. To choose the other tuning parameter $\lambda$, we employ the five-fold cross validation, where the smoothed empirical loss function $\hat{R}_h$ defined in Equation \eqref{eq-smooth}  is used to compare the performance of the estimators generated from different $\lambda$'s. The optimal $\lambda$ that minimizes this empirical loss turns out to be around $1 \times 10^{-6}$. 
	
	We consider three different quantile levels ($\tau = 0.25, 0.5, 0.75$) to compare these two approaches. In particular, for each $\tau$, we define the mean squared integration error of the estimated slope function $\hat{\beta}(\cdot, \tau)$ as $\text{MISE}(\hat{\beta}(\cdot,\tau)) = \int_0^1 \{\hat{\beta}(t, \tau) - \beta(t)\}^2 dt$. 
	Table \ref{tab:MISE} summarizes MISEs of the estimated coefficient functions using these two methods. Our proposed estimator compares favorably with the one using the FPCA method regardless of the distribution of the random error and the quantile level.
	This is also demonstrated by Figure \ref{fig:estimate}. Figure \ref{fig:estimate} displays the average of the estimated coefficient functions across the 200 simulation runs for normally distributed errors under SNR = 10. For these three quantile levels ($\tau = 0.25, 0.5, 0.75$), the mean estimate generated from our proposed approach is nearly unbiased, while the counterpart from the FPCA method does not behave well near the boundaries. 
	
	\begin{table}[H]
		\centering
		
		\addtolength{\tabcolsep}{-1pt}    
		\begin{tabular}{cccccccc}
			\hline
			\multirow{2}{*}{Distribution}&
			\multirow{2}{*}{SNR}&
			\multicolumn{2}{c}{$\tau = 0.25$} &
			\multicolumn{2}{c}{$\tau = 0.5$} &
			\multicolumn{2}{c}{$\tau = 0.75$}  \\
			
			&	& FPCA & RKHS  & FPCA & RKHS  & FPCA & RKHS   \\
			\cline{1-8}
			Normal	&10	& 1.60 (1.58)  & 0.55 (0.57) & 1.47 (1.42) & 0.29 (0.30)  & 1.53 (1.60)  & 0.33 (0.32)    \\
			Normal	& 5	& 2.73 (1.73) & 0.95 (0.98)  & 2.57 (1.44) & 0.53 (0.56)  & 2.85 (2.08) & 0.64 (0.73)   \\
			$t_3$	& 10	& 0.88 (1.05) & 0.37 (0.34)  & 0.61 (0.72) & 0.19 (0.20)  & 0.81 (0.93) & 0.15 (0.16)   \\
			$t_3$ & 5	& 1.72 (1.53)  & 0.58 (0.56)  & 1.39 (1.40) & 0.31 (0.35)  & 1.82 (1.52) & 0.28 (0.32)   \\
			\hline
			
		\end{tabular}
		
		
		\caption{
			Summary of averages and standard errors of MISE ($ \times 10^{-2})$ for the estimated slope functions based on the FPCA method and the RKHS method in different scenarios across 200 simulation runs.}	
		
		\label{tab:MISE}
	\end{table}

	\begin{table}[H]
		\centering
		
		\addtolength{\tabcolsep}{3pt}    
		\begin{tabular}{cccccccc}
			\hline
			\multirow{2}{*}{$\tau$}&
			\multicolumn{3}{c}{SNR = 10} & &
			\multicolumn{3}{c}{SNR = 5}   \\
			\cline{2-4} \cline{6-8}
			& 0.1 & 0.5  & 0.9 & & 0.1  & 0.5 & 0.9   \\
			\cline{1-8}
			0.25	& 0.947 & 0.967 & 0.950 &  & 0.950 & 0.975 & 0.950   \\
			0.5	& 0.957 & 0.947 & 0.930  &  & 0.960 & 0.940 & 0.940    \\
			0.75	& 0.940 & 0.923 & 0.943 &  & 0.925 & 0.925 & 0.945  \\
			
			\hline
			
		\end{tabular}
		
		
		\caption{
			Coverage probabilities of $\beta_0(t, \tau)$ at $t = 0.1, 0.5$ and 0.9 and $\tau \in \{0.25, 0.5, 0.75\}$.}	
		
		\label{tab:CP}
	\end{table}

	\subsection{Confidence intervals/bands for $\beta(\cdot, \tau)$}
\begin{figure}[H]
		\centering
		{\includegraphics[width=14cm]{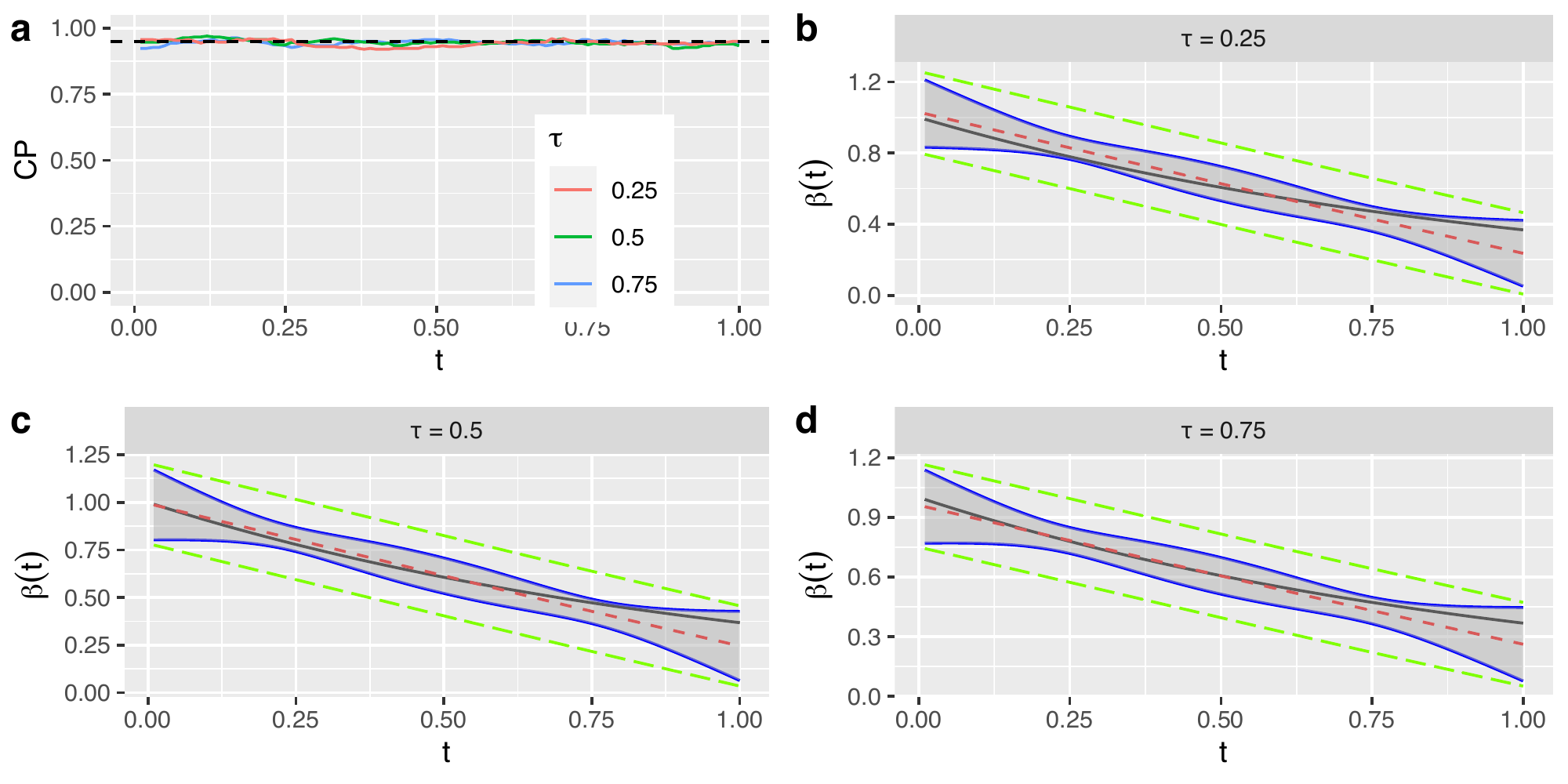}}
		\caption{(a): The true coverage probability of the pointwise confidence intervals across 200 simulations for normally distributed errors with SNR = 10. (b)-(d): In each panel, the solid black line represents the true slope function, while the dashed red line represents the estimated slope function from one random simulation run. The solid blue lines and the dashed green lines represent 95\% pointwise confidence intervals and a 95\% SCB, respectively.}
		\label{fig:confidence}
	\end{figure}

	In this section we study the performance of the pointwise confidence intervals introduced in Remark \ref{re-pointCI} and the SCB in Remark \ref{re-SCB}. 
	Panel (a) in Figure \ref{fig:confidence} displays the coverage probabilities of the pointwise confidence intervals for normally distributed errors with SNR = 10. For all of these three quantile levels, the true coverage probabilities are always close to the nominal level 95\%. Panels (b) - (d) depict the estimated coefficient function based on our proposed method and the corresponding pointwise confidence intervals and the SCB from one random instance for the three quantile levels. The estimated coefficient function provides a reasonable estimate of the true coefficient function for all three quantile levels and the SCB fully covers the true slope function. 
	Table \ref{tab:CP} details the coverage probability of $\hat{\beta}(t, \tau)$ for $\beta_0(t, \tau)$ at $t = \{0.1, 0.5, 0.9\}$ and $\tau = \{0.25, 0.5, 0.75\}$. Even with a moderate number of training curves ($n = 200$), the true coverage probabilities are close to the nominal level 95\% for both of these two signal to noise levels.

	\subsection{Inference on the conditional quantile}
	We investigate the confidence intervals for conditional quantile proposed in Section \ref{sec-CIquant}. In addition to the training sample that consists of $n = 200$ or 400 curves and scalar responses, we generate a new observation $x_0(t)$ and its corresponding response variable $y_0$; this new observation pair is fixed across $B = 200$ simulation runs. To evaluate the converage probability of the constructed confidence interval for the conditional quantile given $x_0$, in each simulation run we implement the proposed estimation procedure on the independently generated training sample to generate a 95\% confidence interval. The true coverage probability is approximated by the frequency of covering the true conditional quantile of $y_0$ given $x_0$ across 200 simulations. Table \ref{tab:CP-quant} summarizes the coverage probabilities of the constructed confidence interval in different designs. Regardless of sample size or SNR levels, the true coverage probability is always close to the nominal level 95\%. When the SNR is fixed, increasing the size of the training sample can yield confidence intervals for the conditional quantile with a true coverage probability closer to the nominal level.

	\begin{table}[H]
		\centering
		
		\addtolength{\tabcolsep}{3pt}    
		\begin{tabular}{cccccc}
			\hline
			\multirow{2}{*}{$\tau$}&
			\multicolumn{2}{c}{$n = 200$} & &
			\multicolumn{2}{c}{$n = 400$}   \\
			\cline{2-3} \cline{5-6}
			& SNR = 10 &  SNR = 5 & & SNR = 10 & SNR = 5  \\
			\cline{1-6}
			0.25	& 0.940	&  0.940 &  & 0.950  &  0.950   \\
			0.5	& 0.945	& 0.945  &  & 0.950    & 0.955    \\
			0.75	& 0.930	& 0.920 &  & 0.940  &   0.940  \\
			
			\hline
			
		\end{tabular}
		
		
		\caption{
			Coverage probabilities of conditional quantile.}	
		\label{tab:CP-quant}
	\end{table}

	\section{Real example}
	\label{sec:real}
	In this section,  we apply the proposed method to one real-world example to demonstrate its performance. 
	
	Alzheimer's disease (AD) is the most common cause of dementia, a general term for the impaired cognitive abilities such as memorizing, thinking, or making decisions that interferes with daily life. According to \cite{ferreira2011}, clinicians  usually leverage a history of progressive and characteristic cognitive decline and the presence of objective cognitive deficits in the diagnosis of AD. Examinations leveraging neuroimaging plays a critical role in the diagnostic investigation of dementia. In particular, these examinations can not only facilitate identification of non-AD pathological processes that are related to cognitive decline but also are useful to pinpoint informative biomarkers that can help to enhance the diagnosis of AD. Substantial studies have shown correlations between quantitative measures collected by magnetic resonance brain imaging and AD progression \citep{ledig2018}.

	\begin{figure}[H]
		\centering
		{\includegraphics[width=14cm]{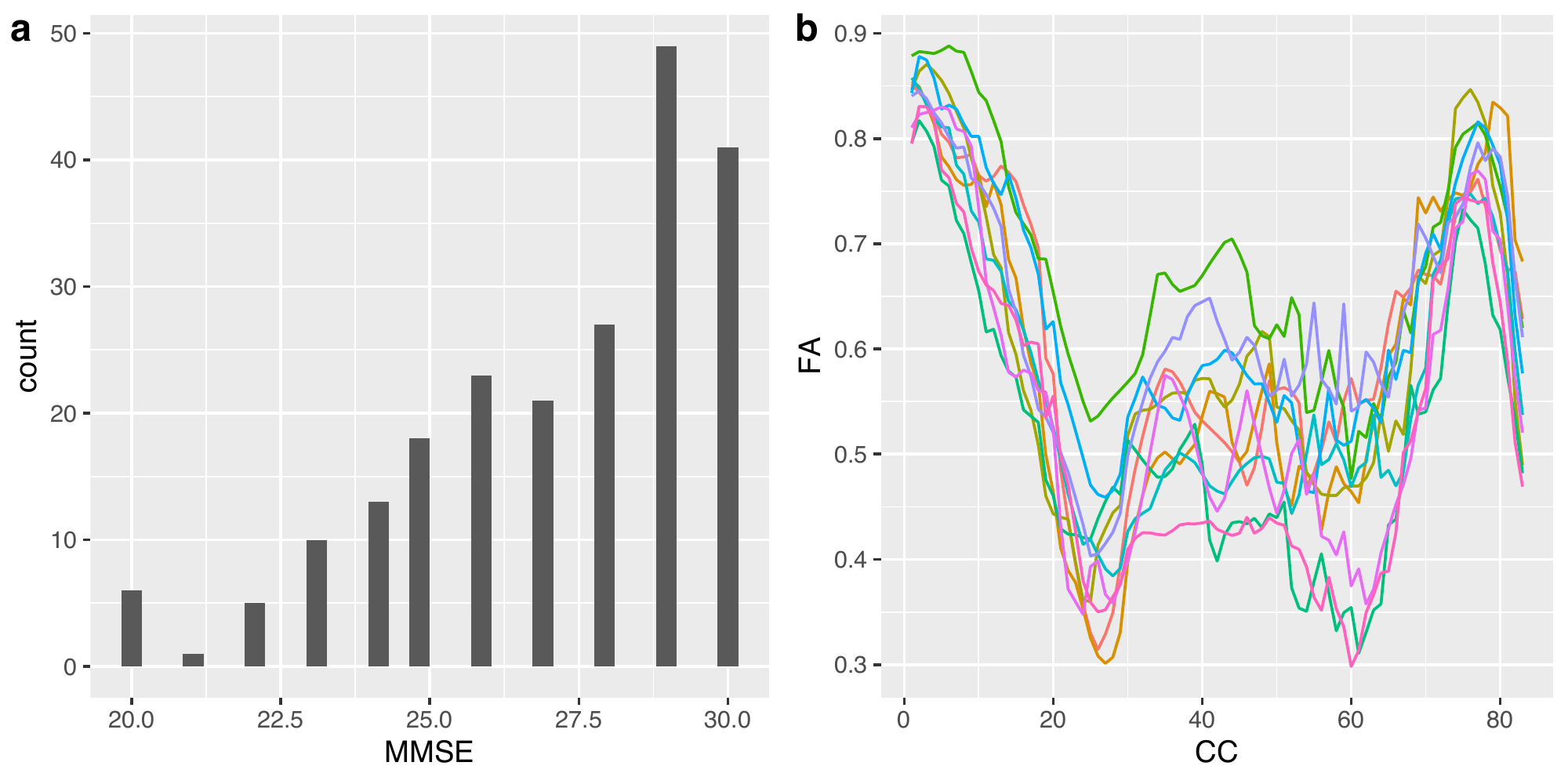}}
		\caption{(a): Histogram of the MMSE scores from $n = 214$ subject. (b): The trajectories of FA values measured at 83 different locations along CC from 10 randomly selected subjects.}
		\label{fig:ADNI}
	\end{figure}
	
	In this study, we obtained the data from the ongoing Alzheimer’s Disease Neuroimaging Initiative (ADNI), that unites researchers working to define the progression of AD. Particularly, the ADNI aims to determine the relationships between clinical, cognitive, imaging, genetic, and biochemical biomarkers over the whole range of AD progression. In this dataset, the functional covariate $X$ is fractional anisotropy (FA) collected by diffusion tensor imaging. The FA values were measured at 83 locations along the corpus callosum (CC) fiber tract for each subject. Panel (b) in Figure \ref{fig:ADNI} depicts the FA trajectories of 10 randomly selected subjects. The scalar response of interest, $Y$, is the mini-mental state examination (MMSE) scores. We chose MMSE as the response since it is one
	of the most widely used test of cognitive functions for assessing the level of dementia a patient may have. Various functional regression models have been proposed to model the relations between the MMSE and FA trajectories; see \cite{zhu2012} and \cite{tang2021partial} for instance. 
	Panel (a) in Figure \ref{fig:ADNI} displays the histogram of the MMSE. Obviously the distribution of $Y$ is highly left skewed, hence a simple summary statistic like sample mean cannot adequately describe the (conditional) distribution of $Y$ given $X$. We are concerned about using functional linear quantile regression to provide a comprehensive characterization of the conditional distribution of the MMSE in this study.

	After removing 3 subjects with missing values in FA trajectories, $n = 214$ subjects were retained in our analysis. To compare the prediction performance of
	the FPCA-based method by \cite{kato2012} and our RKHS-based method, a training set with 80\% of all subjects were randomly selected and the remaining subjects constituted the test set. To better assess prediction accuracy, we randomly splitted the whole dataset 500 times with this ratio. 
	Table \ref{tab:ADNI} summarizes the mean prediction error, which is the average loss function ($\rho_{\tau}(\cdot)$) calculated on the test set, and the standard errors across the 500 splits for these two methods. Our proposed method outperforms its counterpart in terms of prediction accuracy in all of these three quantile levels.

	\begin{table}[H]
		\centering
		
		\begin{tabular}{cccc}
			\hline
			Method&{$\tau = 0.25$} &{$\tau = 0.5$} & {$\tau = 0.75$}  \\
			
			\hline
			
			FPCA	& 2.63 (1.46)	& 2.41 (0.95)  & 0.80 (0.25)   \\
			RKHS	& 0.90 (0.11)	& 1.01 (0.11) & 0.69 (0.08) \\
			\hline
			
		\end{tabular}
		\caption{
			Summary of averages and standard errors of prediction errors based on the FPCA method and the RKHS method across the 500 random splits.}	
		
		\label{tab:ADNI}
	\end{table}
	
\begin{figure}[H]
		\centering
		{\includegraphics[width=14cm]{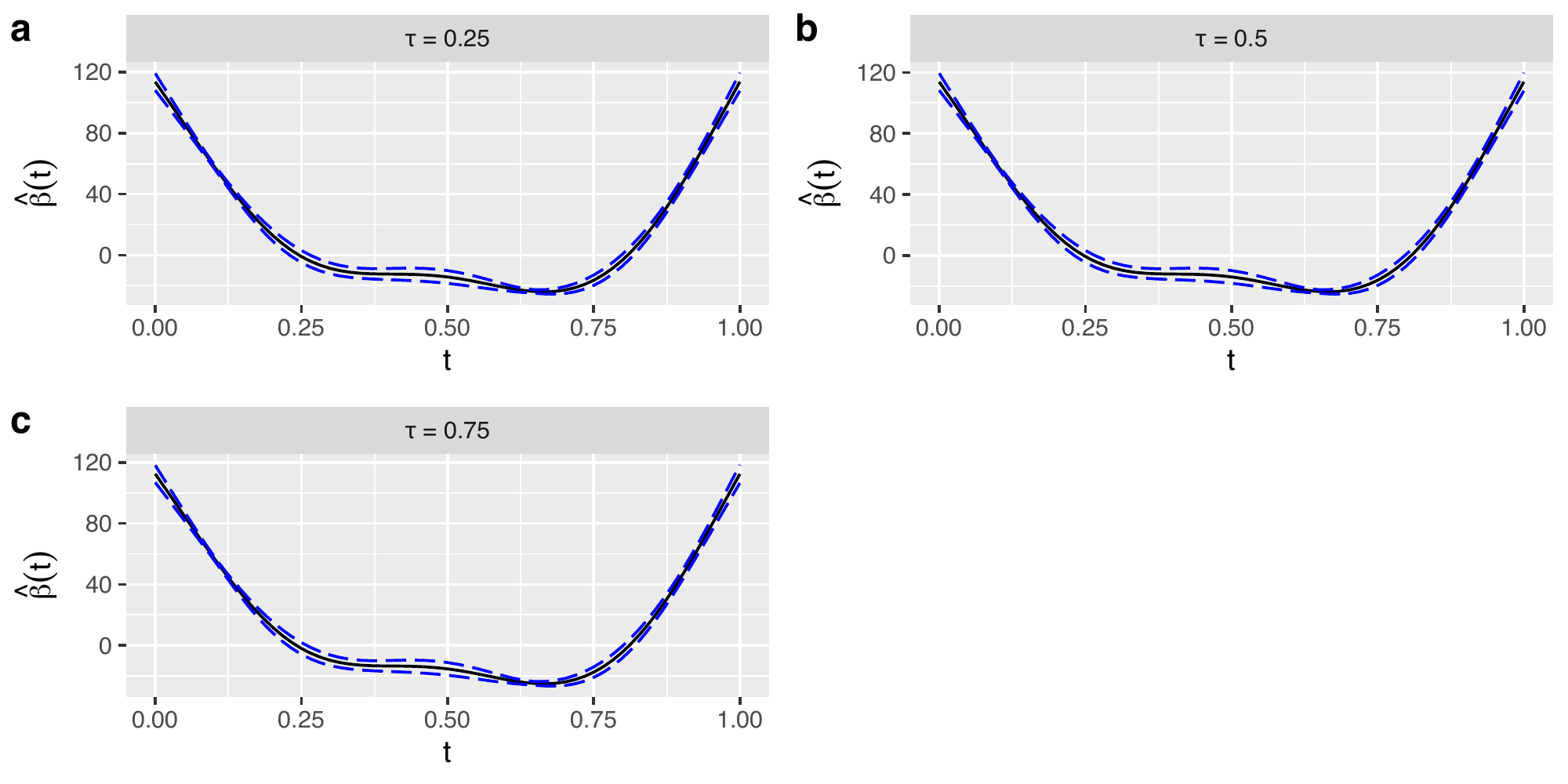}}
		\caption{Estimated coefficient functions and the corresponding 95\% pointwise confidence intervals for $\tau \in \{0.25, 0.5, 0.75\}$.}
		\label{fig:ADNI_beta}
	\end{figure}

	Figure \ref{fig:ADNI_beta} depicts the estimated slope functions from minimizing Equation \eqref{eq-repre} for the above three quantile levels. In addition, the corresponding 95\% pointwise confidence intervals for each estimated slope function are also displayed in Figure \ref{fig:ADNI_beta}. Note $t$ denotes the rescaled CC in Figure \ref{fig:ADNI}. 
	We follow the idea of \cite{yao2005b} to employ
	the one-curve-leave-out analysis to obtain pointwise confidence intervals, which are depicted in Figure \ref{fig:ADNI_beta}. These intervals
	suggest that the FA trajectory has a time-vary effect on the conditional quantile of the MMSE scores, and this effect becomes weaker at the middle range of the CC fiber tract. 
	
	\section{Conclusions}
	\label{sec:conclusion}
	In this paper we employ the kernel convolution technique to smooth the ordinary loss function in FLQR. Under the framework of RKHS, we establish a functional Bahadur representation for the minimizer of the regularized and smoothed empirical risk function. Then we develop inference tools such as pointwise confidence intervals and a SCB for the coefficient function. With the aid of the representer theorem, we design a gradient descent algorithm to facilitate estimation of the intercept and the coefficient function in FLQR. This algorithm is essentially different from that based on FPCA proposed by \cite{kato2012}. Furthermore, empirical studies demonstrate that our proposed estimator outperforms the latter in finite sample.

	Traditional theoretical work on FLQR was focused on estimation consistency and/or convergence rates. For instance, \cite{kato2012} studied the minimax rate in the estimation of the coefficient function and the conditional quantile function. \cite{yao2017} established the convergence rate of the estimated coefficient function in a partially FLQR model where a functional covariate as well as a finite-dimensional covariate is considered. In contrast, little work has been done to develop statistical inferences for this model. 
	Though \cite{li2016inference} is an exception, this work imposes strong assumptions on the coefficient function. In contrast, mild conditions are required for it under the theoretical framework developed in this paper.

	\setcounter{equation}{0}
	\setcounter{figure}{0}
	\setcounter{table}{0}
	\setcounter{page}{1}
	\makeatletter
	\renewcommand{\theequation}{A\arabic{equation}}
	\renewcommand{\thefigure}{A\arabic{figure}}
	\renewcommand{\bibnumfmt}[1]{[A#1]}
	\renewcommand{\citenumfont}[1]{A#1}
	
	\setcounter{section}{0}
	\renewcommand{\thesection}{A.\arabic{section}}
	
\appendix	
\section{Proof of Lemmas and Theorems}
The section contains the proofs of the lemma and theorems in the main manuscript. 
	
\subsection{Proof of Lemma \ref{lemma1}}

	Under Assumption  \ref{ass-f}, we obtain a Taylor expansion
	$$
	f(u + hz | x) = \sum_{j = 0}^s f^{(j)} (u | x) \frac{(hz)^j}{j!} + \frac{(hz)^s}{(s - 1)!} \int_0^1 (1 - w)^{s-1}\left[f^{(s)}(u + whz | x) - f^{(s)}(u | x)\right] dw.
	$$
	Recall that for any $x \in \mbox{supp}(X)$, 
	\begin{align}
		\nonumber
		\E[k_h(u - Y) | x] - f(u | x) & = \int k_h(u - y) f(y | x) dy - f(u | x) = \int k(z) \left[f(u + hz | x) - f(u | x) \right] dz \\
		\label{eq-remainder}
		& = \int_0^1 (1 - w)^{s - 1} \int \frac{(hz)^s}{(s - 1)!} k(z) \left[f^{(s)}(u + whz | x) - f^{(s)}(u | x)\right] dz dw 
	\end{align}
	(i)	As the check function can be written as
	$$
	\int \rho_{\tau}(u) dG(u) = (1 - \tau) \int_{-\infty}^0 G(v)dv + \tau \int_0^{\infty} [1 - G(u)]du
	$$
	
	for any cdf $G$. Therefore,
	$$
	R(\theta; \tau) = \int \left\{(1 - \tau) \int_{-\infty}^0 \int_{-\infty}^{v + \inner{R \lo {x}}{ \theta}} f(u | x)du dv + \tau  \int^{\infty}_0 \int^{\infty}_{v + \inner{R \lo {x}}{ \theta}} f(u | x)du dv\right\} d P_X(x).
	$$
	Similarly, we have
	\begin{align*}
		R_h(\theta; \tau) & = 
		\int \left\{(1 - \tau) \int_{-\infty}^0 \int_{-\infty}^{v + \inner{R \lo {x}}{ \theta}} \E[k_h(u - Y) | x] du dv \right.\\
		&\quad\quad +  \left.\tau  \int^{\infty}_0 \int^{\infty}_{v + \inner{R \lo {x}}{ \theta}} \E[k_h(u - Y) | x] du dv\right\} d P_X(x).
	\end{align*}
	
	Under Assumption \ref{ass-K}, integrating \eqref{eq-remainder} yields
	\begin{align*}
		L_1 & := \left| \int_{-\infty}^0 \int_{-\infty}^{v + \inner{R \lo {x}}{ \theta}} \left\{\E[k_h(u - Y) | x] - f(u | x)  \right\} dudv \right| \\
		& = \left|\int_0^1 (1-w)^{s-1} \int \frac{(hz)^s}{(s - 1)!} k(z) \int_{-\infty}^0 \int_{-\infty}^{v + \inner{R \lo {x}}{ \theta}}\left[f^{(s)}(u + whz | x) - f^{(s)}(u | x)\right] du dv dz dw \right|\\
		& = \left| \int_0^1 (1-w)^{s-1} \int \frac{(hz)^s}{(s - 1)!} k(z)   \left[f^{(s-2)}(\inner{R \lo {x}}{ \theta} + whz | x) - f^{(s-2)}(\inner{R \lo {x}}{ \theta} | x)\right]  dz dw\right| \\
		& \leq Ch^{s+1}
	\end{align*}
	since $\int |z^{s+1}k(z)| dz < \infty$ by Assumption \ref{ass-K} and that $f^{(s-2)}(\cdot | x)$ is Lipschitz by Assumption \ref{ass-f}. Using the same technique, we are able to show that 
	$$
	\left| \int^{\infty}_0 \int^{\infty}_{v + \inner{R \lo {x}}{ \theta}} \left\{\E[k_h(u - Y) | x] - f(u | x)  \right\} dudv \right| \leq Ch^{s+1}
	$$ for some constant $C$. This establishes (i). \\

	\noindent (ii) By the Lebesgue dominated convergence theorem, we have
	$$
	R^{(1)}(\theta; \tau)\theta_1 = \E \left[\big(F(\inner{R \lo {X}}{ \theta} | X) - \tau\big) \inner{R \lo {X}}{ \theta_1} \right] = \int \left[\int_{-\infty}^{\inner{R \lo {x}}{ \theta}} f(y | x)dy - \tau\right] \inner{R \lo {x}}{ \theta_1} d P_X(x),
	$$
	and 
	\begin{align*}
		R^{(1)}_h(\theta; \tau)\theta_1 & = \E \left\{\left[\bK\left(\frac{\inner{R \lo {X}}{ \theta} - Y}{h}\right) - \tau\right] \inner{R \lo {X}}{ \theta_1}\right\}  \\
		&	= \int \left\{\int_{-\infty}^{\inner{R \lo {x}}{ \theta}} \E[k_h(u - Y | x)]du - \tau\right\} \inner{R \lo {x}}{ \theta_1} d P_X(x).
	\end{align*}
	As $\int u^s k(u)du = 0$ and $\int |u^{s+1}k(u)| du < \infty$, by integrating \eqref{eq-remainder} we have
	\begin{align*}
		L_2 & : = \bigg|\int_{-\infty}^{\inner{R \lo {x}}{ \theta}}  \E[k_h(u - Y | x)] - f(u | x) du \bigg| \\
		& = \left| \int_0^1 (1-w)^{s-1} \int \frac{(hz)^s}{(s - 1)!} k(z)\int_{-\infty}^{\inner{R \lo {x}}{ \theta}}\left[f^{(s)}(u + whz | x) - f^{(s)}(u | x)\right]  dudz dw\right| \\
		& = \left| \int_0^1 (1-w)^{s-1} \int \frac{(hz)^s}{(s - 1)!} k(z)\left[f^{(s-1)}(\inner{R \lo {x}}{ \theta} + whz | x) - f^{(s-1)}(\inner{R \lo {x}}{ \theta} | x)\right]  dz dw\right| \\
		& = \left| \int_0^1 (1-w)^{s-1} \int \frac{(hz)^{s+1}}{(s - 1)!} k(z)\int_0^1 f^{(s)}(\inner{R \lo {x}}{ \theta} + twhz | x) dt  dz dw\right| \\
		& \leq Ch^{s+1},
	\end{align*}
	uniformly given the fact that $f^{(s)}$ is bounded. 
	
	It follows that 
	\begin{align*}
		R^{(1)}_h(\theta; \tau)\theta_1 - R^{(1)}(\theta; \tau)\theta_1 & \leq Ch^{s+1} \E[|\inner{R \lo {X}}{ \theta_1}|]  \\
		& \leq  Ch^{s+1} \left[\E\left(\langle R_X, \theta_1 \rangle^4\right)\right]^{1/4} \leq C h^{s+1} \|\theta_1\|
	\end{align*}
	by the Cauchy-Schwarz inequality and Lemma S.3 of \cite{shang2015b}. This establish (ii). \\
	
	\noindent (iii) The second order Fr\'{e}chet derivative operator of $R(\theta; \tau)$ satisfies
	$$
	R^{(2)}(\theta; \tau)\theta_1\theta_2 = \E \left[f(\inner{R \lo {X}}{ \theta} | X) \inner{R \lo {X}}{ \theta_1} \inner{R \lo {X}}{ \theta_2}\right] = \int f(\inner{R \lo {x}}{ \theta} | x )  \inner{R \lo {x}}{ \theta_1} \inner{R \lo {x}}{ \theta_2} d P_X(x),
	$$
	and, analogously, 
	\begin{align*}
		R^{(2)}_h(\theta; \tau)\theta_1\theta_2 & = \E \left[k_h(\inner{R \lo {X}}{ \theta} - Y) \inner{R \lo {X}}{ \theta_1} \inner{R \lo {X}}{ \theta_2}\right] \\
		& = \int \E\left[k_h(\inner{R \lo {x}}{ \theta} - Y | x )\right]  \inner{R \lo {x}}{ \theta_1} \inner{R \lo {x}}{ \theta_2} d P_X(x).
	\end{align*}
	Setting $u = \inner{R \lo {x}}{ \theta}$. In light of \eqref{eq-remainder}, we have 
	\begin{align*}
		& \left|\E\left[k_h(u - Y) | x\right] - f(u | x) \right| \\
		& \leq Ch^s \int |z^sk(z)| \sup_{(y, x) \in \Re \times \mbox{supp}(X)} \sup_{t: |t| \leq hz} \left|f^{(s)}(y + t | x) - f^{(s)}(y | x)\right|dz \\
		& = o(h^s)
	\end{align*}
	under Assumption \ref{ass-f} and Assumption \ref{ass-K}. 
	
	Then following the proof of part (ii), it is easy to verify that \\
	$|R^{(1)}_h(\theta; \tau)\theta_1\theta_2 - R^{(1)}(\theta; \tau)\theta_1\theta_2| = o(h^s)\|\theta_1\|\|\theta_2\|$. This establishes (iii).  \\
	
	\noindent (iv) Note that 
	\begin{align*}
		R^{(2)}_h(\theta; \tau)\theta_1\theta_2 & = \E \left[k_h(\inner{R \lo {X}}{ \theta} - Y) \inner{R \lo {X}}{ \theta_1} \inner{R \lo {X}}{ \theta_2}\right] \\
		&	= \int k(z) \int f(\inner{R \lo {x}}{ \theta} + hz | x) \inner{R \lo {x}}{ \theta_1} \inner{R \lo {x}}{ \theta_2} d P_X(x) dz. 
	\end{align*}
	Since $f(\cdot | x)$ is Lipschitz under Assumption \ref{ass-f} (d), 
	\begin{align*}
		& \left| R_h^{(2)}(\theta + \theta^*; \tau)\theta_1\theta_2 - R_h^{(2)}(\theta; \tau)\theta_1\theta_2  \right| \\
		& \leq C \int |k(z)| \int |\inner{R \lo {X}}{ \theta^*}| |\inner{R \lo {X}}{ \theta_1}| |\inner{R \lo {X}}{ \theta_2}| d P_X(x) dz\\
		& \leq C \E[|\inner{R \lo {X}}{ \theta^*}| |\inner{R \lo {X}}{ \theta_1}| |\inner{R \lo {X}}{ \theta_2}|] \\
		& \leq C \left[\E\left(\langle R_X, \theta\hi * \rangle^4\right)\right]^{1/4} \left[\E\left(\langle R_X, \theta_1 \rangle^4\right)\right]^{1/4} \left[\E\left(\langle R_X,\theta_2\rangle^4\right)\right]^{1/4} \\
		& \leq C \|\theta\hi*\| \|\theta_1\| \|\theta_2\|
	\end{align*}
	uniformly in $(\theta, h, \theta^*, \tau)$. 
	
	\subsection{Proof of Theorem \ref{th-thlam}}

	Recall that $R(\theta; \tau) = \int \rho_{\tau}(u) \E[d \hat{F}(u; \theta)]$, where the expectation is taken with respect to the joint distribution of $(X, Y)$. 
	Since $\hat{F}(u; \theta) = \frac{1}{n} \sum_{i = 1}^n \mathbbm{1}(\eps_i(\theta) \leq u)$, $\E \hat{F}(u; \theta) = P(\eps(\theta) \leq u)$, where $e(\theta) = Y - \inner{R_X}{\theta}$. Hence $\frac{d \E[\hat{F}(u; \theta)]}{du} \equiv f_{\epsilon}(u; \theta)$, where $f_{\epsilon}(\cdot, \theta)$ denotes the marginal density of $\epsilon(\theta): = Y - \inner{R_X}{\theta}$. It follows that $R(\theta; \tau) = \int \rho_{\tau}(u) f_{\eps}(u; \theta)$. Similarly, $\hat{R}_h(\theta; \tau) = \frac{1}{n} \sum_{i = 1}^n \int \rho_{\tau}(u) k_h(u - \eps_i(\theta)) du$ and $R_h(\theta; \tau) = \E[\hat{R}_h(\theta; \tau)] = \int \rho_{\tau}(u) \E[k_h(u - \eps(\theta))] du$.
	Let $\ti{K}_h(u; \theta) = \E[k_h(u - \eps(\theta))]$, where the expectation is taken with respect to the distribution of $\eps(\theta)$. In other words, $\ti{K}_h(u; \theta)
	= \int f_{\eps}(v; \theta) k_h(u - v)dv$ and $R_h(\theta; \tau) = \int \rho_{\tau}(u)\ti{K}_h(u; \theta) du$. So $\ell_{h\lambda}(\theta) = \int \rho_{\tau}(u)\ti{K}_h(u; \theta) du + \frac{\lambda}{2} J(\beta, \beta)$ based on the definition of $\ell_{h\lambda}$ in \eqref{eq-PPL}. 
	
	The first-order Fr\'{e}chet derivative of $\ell_{h\lambda}$ w.r.t $\theta$ satisfies that for any $\Delta\theta = (\Delta\alpha, \Delta\beta) \in \ca H$, 
	$$
	S_{h\lambda}(\theta) \Delta\theta \equiv  D_{\theta} \ell_{h\lambda}(\theta) \Delta\theta = \int \rho_{\tau}(u)  D_{\theta} \ti{K}_h(u; \theta)\Delta\theta  du +  \inner{P_{\lambda}\theta}{\Delta\theta},
	$$
	where $P_{\lambda}(\theta) = (0, W_{\lambda}\beta) \in \ca H$ and satisfies that $\inner{P_{\lambda}(\theta)}{\ti{\theta}} = \inner{W_{\lambda}\beta}{\ti{\beta}}_1$ for any $\ti{\theta} = (\ti\alpha, \ti\beta) \in \ca H$ from Proposition 2.5 of \cite{shang2015}. The second-order  Fr\'{e}chet derivative of $\ell_{h\lambda}$ is
	$$
	D^2_{\theta} \ell_{h\lambda}(\theta) \Delta\theta_1\Delta\theta_2 = \int \rho_{\tau}(u)  D_{\theta}^2 \ti{K}_h(u; \theta)\Delta\theta_1\Delta\theta_2  du + \inner{P_{\lambda}\Delta\theta_1}{\Delta\theta_2}. 
	$$
	For notational brevity, we shall drop $\tau$ from $\theta_{0h}(\tau)$ if the quantile level is clear from the context. 
	Since $\theta_{0h}$ minimizes $R_h(\theta; \tau)$, $S_{h\lambda}(\theta_{0h}) = P_{\lambda}\theta_{0h}$. 
	
	Define $T_1(\theta) = \theta - S_{h\lambda}(\theta_0 + \theta)$ for $\theta \in \ca H$, and let $S_{0\lambda}(\theta_0) = P_{\lambda}\theta_0$. Write $T_1(\theta)$ as $T_1(\theta) + S_{0\lambda}(\theta_0) - S_{0\lambda}(\theta_0)$. So $\|T_1(\theta)\| \leq \|T_1(\theta) + S_{0\lambda}(\theta_0) \| + \| S_{0\lambda}(\theta_0) \|$. For the first term, we have
	\begin{align*}
		T_1(\theta) +  S_{0\lambda}(\theta_0) & = \theta - S_{h\lambda}(\theta_0 + \theta) + S_{0\lambda}(\theta_0) \\
		& = \theta - \int \rho_{\tau}(u)  D_{\theta} \ti{K}_h(u; \theta + \theta_0) du - P_{\lambda}(\theta + \theta_0) \\
		&~~~~ + \int \rho_{\tau}(u)D_{\theta} f_{\eps}(u; \theta_0)du + P_{\lambda}\theta_0 \\
		& = \theta - \int \rho_{\tau}(u) [D_{\theta} \ti{K}_h(u; \theta + \theta_0)  - D_{\theta} \ti{K}_h(u; \theta_0) ]du  \\
		&~~~~ - \int \rho_{\tau}(u) D_{\theta} \ti{K}_h(u; \theta_0) du - P_{\lambda}\theta + \int \rho_{\tau}(u)D_{\theta} f_{\eps}(u; \theta_0)du \\
		& = \theta - \int \rho_{\tau}(u)  \left[\int_0^1 D_{\theta}^2 \ti{K}_h(u; s\theta + \theta_0) \theta ds  \right] du \\
		&~~~~ - \int \rho_{\tau}(u) [D_{\theta} \ti{K}_h(u; \theta_0) - D_{\theta} f_{\eps}(u; \theta_0)] du - P_{\lambda}\theta. 
	\end{align*}
	By the proof of Proposition 3.5 of \cite{shang2015}, we know that \\$\left[\int \rho_{\tau}(u)D^2_{\theta}f_{\eps}(u; \theta_0)du + P_{\lambda}\right]\theta = \theta$.
	So 
	\begin{align*}
		T_1(\theta) + S_{0\lambda}(\theta_0) & = \int \rho_{\tau}(u)D^2_{\theta}f_{\eps}(u; \theta_0)\theta du + P_{\lambda}\theta - \int \rho_{\tau}(u)  \left[\int_0^1 D_{\theta}^2 \ti{K}_h(u; s\theta + \theta_0) \theta ds \right] du \\
		&~~~~ - \int \rho_{\tau}(u) [D_{\theta} \ti{K}_h(u; \theta_0) - D_{\theta} f_{\eps}(u; \theta_0)] du - P_{\lambda}\theta \\
		& = \int \rho_{\tau}(u) \left(\int_0^1[D^2_{\theta}f_{\eps}(u; \theta_0)\theta - D_{\theta}^2 \ti{K}_h(u; s\theta + \theta_0) \theta]ds\right)du \\
		&~~~~ - \int \rho_{\tau}(u) [D_{\theta} \ti{K}_h(u; \theta_0) - D_{\theta} f_{\eps}(u; \theta_0)] du =: \RN{1} + \RN{2}.
	\end{align*}
	For $\RN{1}$, we have
	\begin{align*}
		& \norm{\int \rho_{\tau}(u) [D_{\theta}^2 \ti{K}_h(u; s\theta + \theta_0) - D^2_{\theta}f_{\eps}(u; \theta_0)] du } _{\mbox{op}} \\
		&~~= \norm{\int \rho_{\tau}(u) [D_{\theta}^2 \ti{K}_h(u; s\theta + \theta_0)  - D^2_{\theta}f_{\eps}(u; s\theta + \theta_0) + D^2_{\theta}f_{\eps}(u; s\theta + \theta_0) - D^2_{\theta}f_{\eps}(u; \theta_0)]}_{\mbox{op}} \\
		&~~\leq \norm{\int \rho_{\tau}(u) [D_{\theta}^2 \ti{K}_h(u; s\theta + \theta_0)  - D^2_{\theta}f_{\eps}(u; s\theta + \theta_0)]du}_{\mbox{op}}  \\
		&\quad\quad\quad + \norm{\int \rho_{\tau}(u) [D^2_{\theta}f_{\eps}(u; s\theta + \theta_0) - D^2_{\theta}f_{\eps}(u; \theta_0)]du}_{\mbox{op}} 
	\end{align*}
	Note that
	\begin{align*}
		& \norm{\int \rho_{\tau}(u) [D_{\theta}^2 \ti{K}_h(u; s\theta + \theta_0)  - D^2_{\theta}f_{\eps}(u; s\theta + \theta_0)]du}_{\mbox{op}} \\
		& = \sup_{\|\Delta\theta_1\| = \|\Delta\theta_2\|=1} \bigg|\int \rho_{\tau}(u)D_{\theta}^2 \ti{K}_h(u; s\theta + \theta_0) \Delta\theta_1\Delta\theta_2du - \int \rho_{\tau}(u)D^2_{\theta}f_{\eps}(u; s\theta + \theta_0)\Delta\theta_1\Delta\theta_2du\bigg|. 
	\end{align*}
	For any $\Delta\theta_1, \Delta\theta_2$ with $\|\Delta\theta_1\| = \|\Delta\theta_2\|=1$, based on the proof of (iii) of Lemma \ref{lemma1}, we have
	\begin{align*}
		&\bigg|\int \rho_{\tau}(u)D_{\theta}^2 \ti{K}_h(u; s\theta + \theta_0) \Delta\theta_1\Delta\theta_2du - \int \rho_{\tau}(u)D^2_{\theta}f_{\eps}(u; s\theta + \theta_0)\Delta\theta_1\Delta\theta_2du\bigg|\\
		&=\bigg|\int \E\left[k_h(\inner{R \lo {x}}{ s\theta + \theta_0} - Y | x )\right]  \inner{R \lo {x}}{ \Delta\theta_1} \inner{R \lo {x}}{ \Delta\theta_2} d P_X(x) \\
		& ~~ - \int f(\inner{R \lo {x}}{ s\theta + \theta_0} | x )  \inner{R \lo {x}}{ \Delta\theta_1} \inner{R \lo {x}}{ \Delta\theta_2} d P_X(x)\bigg|\\
		& = o(h^s) \E \left(| \inner{R \lo {X}}{ \Delta\theta_1} \inner{R \lo {X}}{ \Delta\theta_2}|\right) \leq o(h^s)  \left[\E_X(|\inner{R \lo {X}}{ \Delta\theta_1}|^2)\right]^{1/2}\left[\E_X(|\inner{R \lo {X}}{ \Delta\theta_1}|^2)\right]^{1/2}
	\end{align*}
	By Lemma S.3 of \cite{shang2015b},the last term is bounded from above by a constant if \eqref{eq-moment} in Assumption \ref{ass-X} is met. Similarly, for any $\Delta\theta_1, \Delta\theta_2$ with $\|\Delta\theta_1\| = \|\Delta\theta_2\|=1$
	\begin{align*}
		&\bigg| \int \rho_{\tau}(u)D^2_{\theta}f_{\eps}(u; s\theta + \theta_0)\Delta\theta_1\Delta\theta_2du - \int \rho_{\tau}(u)D^2_{\theta}f_{\eps}(u; \theta_0)\Delta\theta_1\Delta\theta_2du \bigg|\\
		&=\bigg|\int f(\inner{R \lo {x}}{ s\theta + \theta_0} | x )  \inner{R \lo {x}}{ \Delta\theta_1} \inner{R \lo {x}}{ \Delta\theta_2} d P_X(x)  \\
		& \quad \quad - \int f(\inner{R \lo {x}}{\theta_0} | x )  \inner{R \lo {x}}{ \Delta\theta_1} \inner{R \lo {x}}{ \Delta\theta_2} d P_X(x)\bigg|\\
		& \leq C_3\E \left(|\inner{R \lo {X}}{\theta} \inner{R \lo {X}}{ \Delta\theta_1} \inner{R \lo {X}}{ \Delta\theta_2}|\right) \leq C_4\|\theta\|
	\end{align*}
	with some positive constant $C$ independent of $n$. Hence $\|\RN{1}\| \leq C[o(h^s) + \|\theta\|] \cdot \|\theta\|$ uniformly over $\tau \in [\underline{\tau}, \bar{\tau}]$. 
	
	We treat $\RN{2}$ next. For any $\Delta\theta$ with $\|\Delta\theta\| = 1$, by the proof of (ii) of Lemma \ref{lemma1}, we have
	\begin{align*} 
		&\bigg|\int \rho_{\tau}(u)D_{\theta} \ti{K}_h(u; \theta_0) \Delta\theta du - \int \rho_{\tau}(u)D_{\theta}f_{\eps}(u; \theta_0)\Delta\theta du\bigg|\\
		& = \bigg|\int \left\{\int_{\infty}^{\inner{R \lo {x}}{ \theta}} \E[k_h(u - Y | x)]du - \tau\right\} \inner{R \lo {x}}{ \Delta\theta} d P_X(x)  \\
		& ~~~ -  \int \left[\int_{\infty}^{\inner{R \lo {x}}{ \theta}} f(y | x)dy - \tau\right] \inner{R \lo {x}}{ \theta_1} d P_X(x)\bigg| \\
		& \leq Ch^{s+1}  \E[|\inner{R \lo {X}}{\Delta \theta}|]  = O(h^{s+1}). 
	\end{align*}
	It follows that there exists some $C > 0$ such that $\|T_1(\theta) + S_{0\lambda}(\theta_0) \| \leq C[o(h^s) + \|\theta\|] \cdot \|\theta\| + Ch^{s+1}$. 
	
	Lastly, we handle $S_{0\lambda}(\theta_0) = P_{\lambda}\theta_0$. Write $\beta_0 = \sum_{\nu} b_{\nu}\hi 0\phi_{\nu}$. According to Proposition \ref{propW}, $W_{\lambda}\beta_0 = \sum_{\nu} b_{\nu}\hi 0 \frac{\lambda\rho_{\nu}}{1+ \lambda\rho_{\nu}} \phi_{\nu}$. Hence $||S_{0\lambda}(\theta_0) ||^2 = \inner{W_{\lambda}\beta_0} {W_{\lambda}\beta_0} _1 = \sum_{\nu} |b_{\nu}^0|^2 \frac{\lambda^2\rho_{\nu}^2}{(1+ \lambda\rho_{\nu})^2}(1 + \lambda\rho_{\nu}) = \lambda \sum_{\nu} |b_{\nu}^0|^2\rho_{\nu} \frac{\lambda\rho_{\nu}}{1+ \lambda\rho_{\nu}} = o(\lambda)$. The last equation holds since $J(\beta_0, \beta_0) = \sum_{\nu} |b_{\nu}^0|^2\rho_{\nu} < \infty$, and $\sum_{\nu}   |b_{\nu}^0|^2\rho_{\nu} \frac{\lambda\rho_{\nu}}{1+ \lambda\rho_{\nu}}  = o(1)$ as $\lambda \rightarrow 0$ by the dominated convergence theorem. Therefore, $||S_{0\lambda}(\theta_0)|| = o(\lambda^{1/2})$. 
	Let $\mathbb{B}(r) = \{\theta \in \ca H: \|\theta\| \leq r\}$ be the ball in $\ca H$ of radius $r$. If we take $r_{1n} = C\max(\lambda^{1/2}, h^{s+1})$, which is of the order $\lambda^{1/2}$ if $h = O(\lambda^{1/(2s+2)})$. 
	Previous results indicate that $T_1(\mB(r_{1n})) \subset \mB(r_{1n})$. 
	
	Next we show that $T_1$ is a contraction mapping. For any $\theta_j = (\alpha_j, \beta_j) \in \mB(r_{1n})$ for $j = 1, 2$, Taylor's expansion yields that
	\begin{align*}
		&\quad T_1(\theta_1) - T_1(\theta_2) \\
		& = \theta_1 - \theta_2 - \int \rho_{\tau}(u)  D_{\theta} \ti{K}_h(u; \theta_1 + \theta_0) du + \int \rho_{\tau}(u)  D_{\theta} \ti{K}_h(u; \theta_2 + \theta_0) du \\
		& ~~~ - P_{\lambda}(\theta_1 + \theta_0) + P_{\lambda}(\theta_2 + \theta_0)  \\
		& = \theta_1 - \theta_2 - P_{\lambda}(\theta_1 - \theta_2) - \int \rho_{\tau}(u) \left[ \int_0^1 D^2_{\theta}\ti{K}_h(u; \theta_2 + \theta_0 + s(\theta_1 - \theta_2))(\theta_1 - \theta_2)ds\right] du \\
		& = \int \rho_{\tau}(u)D^2_{\theta}f_{\eps}(u; \theta_0)(\theta_1 - \theta_2) du + P_{\lambda}(\theta_1 - \theta_2) -  P_{\lambda}(\theta_1 - \theta_2) \\
		&\quad - \int \rho_{\tau}(u) \left[ \int_0^1 D^2_{\theta}\ti{K}_h(u; \theta_2 + \theta_0 + s(\theta_1 - \theta_2))(\theta_1 - \theta_2)ds\right] du \\
		& = \int \rho_{\tau}(u)  \left\{ \int_0^1 \left[D^2_{\theta}f_{\eps}(u; \theta_0) - D^2_{\theta}\ti{K}_h(u; \theta_2 + \theta_0 + s(\theta_1 - \theta_2))\right](\theta_1 - \theta_2)ds\right\} du. 
	\end{align*}
	By similar arguments as in the above analysis of $T_1(\theta)$, we have
	\begin{align*}
		&\quad \|T_1(\theta_1) - T_1(\theta_2)\| \\
		& \leq \norm{\int \rho_{\tau}(u) \left[\int_0^1D_{\theta}^2 \ti{K}_h(u; \theta_2 + \theta_0 + s(\theta_1 - \theta_2))  - D^2_{\theta}f_{\eps}(u; \theta_2 + \theta_0 + s(\theta_1 - \theta_2))ds\right]du}_{\mbox{op}}  \|\theta_1 - \theta_2\|\\
		& + \norm{\int \rho_{\tau}(u) \left[\int_0^1 D^2_{\theta}f_{\eps}(u; \theta_2 + \theta_0 + s(\theta_1 - \theta_2)) - D^2_{\theta}f_{\eps}(u; \theta_0)ds\right]du}_{\mbox{op}} \|\theta_1 - \theta_2\|\\
		& = o(h^s) \|\theta_1 - \theta_2\| + \norm{\int \rho_{\tau}(u) \left[\int_0^1 D^2_{\theta}f_{\eps}(u; \theta_2 + \theta_0 + s(\theta_1 - \theta_2)) - D^2_{\theta}f_{\eps}(u; \theta_0)ds\right]du}_{\mbox{op}} \|\theta_1 - \theta_2\|. 
	\end{align*}
	Note that for any $\Delta\theta_1, \Delta\theta_2$ with $\|\Delta\theta_1\| = \|\Delta\theta_2\|=1$, 
	\begin{align*}
		& \bigg|\int \rho_{\tau}(u) \left[\int_0^1 D^2_{\theta}f_{\eps}(u; \theta_2 + \theta_0 + s(\theta_1 - \theta_2)) \Delta\theta_1\Delta\theta_2 - D^2_{\theta}f_{\eps}(u; \theta_0) \Delta\theta_1\Delta\theta_2ds\right]du\bigg| \\
		& = \bigg|\int \left[\int_0^1 f(\inner{R \lo {x}}{ \theta_2 + \theta_0 + s(\theta_1 - \theta_2)} | x ) - f(\inner{R \lo {x}}{\theta_0} | x ) ds \right]  \inner{R \lo {x}}{ \Delta\theta_1} \inner{R \lo {x}}{ \Delta\theta_2} d P_X(x)\bigg|\\
		& \leq C_0 \int_0^1 \E \{|\inner{R_X}{\theta_2 + s(\theta_1 - \theta_2}| |\inner{R_X}{\Delta_1}||\inner{R_X}{\Delta_2}|\}ds \\
		& \leq C_0 \int_0^1 [\E \{|\inner{R_X}{\theta_2 + s(\theta_1 - \theta_2}|^4\}]^{1/4}ds \times  [\E \{|\inner{R_X}{\Delta\theta_1}|^4\}]^{1/4}  [\E \{|\inner{R_X}{\Delta\theta_2}|^4\}]^{1/4}\\
		& \leq C_0^{\prime} (\|\theta_2\| + \|\theta_1 - \theta_2\|) < 1/2
	\end{align*}
	because as $n \rightarrow \infty$, $\|\theta_2\| + \|\theta_1 - \theta_2\| = O(r_{1n}) = o(1)$. For the same reason, $o(h^s)\|\theta_1 - \theta_2\| < 1/2 \|\theta_1 - \theta_2\|$. 
	
	Therefore, $T_1$ is contraction mapping on $\mB(r_{1n})$. By contraction mapping theorem, there exists uniquely an element $\theta_{\lambda h}^{\prime} \in \mB(r_{1n})$ such that
	$T_1(\theta_{\lambda h}^{\prime}) = \theta_{\lambda h}^{\prime}$. Define $\theta_{0h,\lambda} = \theta_{\lambda h}^{\prime} + \theta_0$. Then we have $S_{h\lambda}(\theta_{0h,\lambda} ) = 0$ and $\|\theta_{0h,\lambda}  - \theta_0\| \leq r_{1n}$.  This completes the proof.

	\subsection{Proof of Theorem \ref{th-rate}}

	Define the operator $T_2(\theta) = \theta - [DS_{h\lambda}(\theta_{0h,\lambda})]^{-1}S_{n, h\lambda}(\theta_{0h, \lambda} + \theta)$ for $\theta \in \ca H$. We firs that the operator $DS_{h\lambda}(\theta_{0h,\lambda})$ is invertible. For any $\Delta\theta_1, \Delta\theta_2 \in \ca H$, 
	\begin{align*}
		&	|[DS_{h\lambda}(\theta_{0h,\lambda}) - DS_{0\lambda}(\theta_0)]\Delta\theta_1\Delta\theta_2| \\
		& = |[DS_{h\lambda}(\theta_{0h,\lambda}) - DS_{h\lambda}(\theta_{0})]\Delta\theta_1\Delta\theta_2 \\
		&~~+ [DS_{h\lambda}(\theta_{0}) - DS_{0\lambda}(\theta_0)]\Delta\theta_1\Delta\theta_2| \\
		& \leq |R_h^{(2)}(\theta_{0h, \lambda}; \tau)\Delta\theta_1\Delta\theta_2 - R_h^{(2)}(\theta_0; \tau)\Delta\theta_1\Delta\theta_2|  \\
		&\quad + |R_h^{(2)}(\theta_0; \tau)\Delta\theta_1\Delta\theta_2 - R^{(2)}(\theta_0; \tau)\Delta\theta_1\Delta\theta_2| \\
		& \leq C \norm{\theta_{0h, \lambda} - \theta_0} \norm{\Delta\theta_1} \norm{\Delta\theta_2} + o(h^s)  \norm{\Delta\theta_1} \norm{\Delta\theta_2}.
	\end{align*}
	The last equation holds due to (iv) and (iii) of Lemma \ref{lemma1}, respectively. By Theorem \ref{th-thlam}, we have that the operator norm of the difference of $DS_{h\lambda}(\theta_{0h,\lambda})$ and $DS_{0\lambda}(\theta_0)$ is strictly less than 1/2. As shown in the proof of Theorem \ref{th-thlam}, $DS_{0\lambda}(\theta_0) = id$. Therefore 
	$DS_{h\lambda}(\theta_{0h,\lambda})$ is invertible, and the norm of $[DS_{h\lambda}(\theta_{0h,\lambda})]^{-1}$ is between 2/3 and 2. 
	
	Rewrite $T_2$ as
	\begin{align*}
		T_2(\theta) & = -[DS_{h\lambda}(\theta_{0h,\lambda})]^{-1}[DS_{n, h\lambda}(\theta_{0h,\lambda})\theta - DS_{h\lambda}(\theta_{0h,\lambda})\theta] \\
		&\quad\quad -[DS_{h\lambda}(\theta_{0h,\lambda})]^{-1}[S_{n, h\lambda}(\theta_{0h, \lambda} + \theta) - S_{n, h\lambda}(\theta_{0h, \lambda}) - DS_{n, h\lambda}(\theta_{0h, \lambda})\theta] \\
		&\quad\quad -[DS_{h\lambda}(\theta_{0h,\lambda})]^{-1}S_{n,h\lambda}(\theta_{0h,\lambda}) \\
		& : = I_1 + I_2 + I_3. 
	\end{align*}
	
	Next we treat these three terms to verify that $T_2$ is a contraction mapping when its domain is restricted to a specific ball. Let $\bK_h(\cdot)= \bK(\cdot/h)$. 
	To deal with $I_3$, we define $O_i\hi{\prime} = [\bK_h(-\eps_i(\theta_{0h,\lambda})) - \tau]R_{X_i}$ and $O_i = O_i\hi{\prime} - \E(O_i\hi{\prime})$ for $i = 1, \ldots, n$. Since $S_{h\lambda}(\theta_{0h,\lambda}) = 0$, 
	\begin{align*}
		\E[\norm{S_{n,h\lambda}(\theta_{0h,\lambda})}^2] = \E[\norm{S_{n,h\lambda}(\theta_{0h,\lambda}) - S_{h\lambda}(\theta_{0h,\lambda})}^2] = n^{-1} \E[\norm{O_i}^2] \\
		\leq n^{-1} \E\{|\bK_h(-\eps(\theta_{0h,\lambda})) - \tau|^2 \norm{R_X}^2\}. 
	\end{align*}
	To find the term on the right-hand side, we employ the law of iterated expectations. In particular, given $X = x$, 
	\begin{align*}
		\E\{|\bK_h(-\eps(\theta_{0h,\lambda})) - \tau|^2 | X = x\} & = \E[\bK^2_h(\inner{R_X}{\theta_{0h,\lambda}} - Y) | X = x]  \\
		&\quad -2\tau \E[\bK_h(\inner{R_X}{\theta_{0h,\lambda}} - Y) | X = x]  + \tau^2.
	\end{align*}
	Integration by parts leads to 
	\begin{align*}
		& \E[\bK_h(\inner{R_X}{\theta_{0h,\lambda}} - Y) | X = x] \\
		& = \int \bK_h(\inner{R_x}{\theta_{0h,\lambda}} - y) f(y | x) dy \\
		& = \int k_h(\inner{R_x}{\theta_{0h,\lambda}} - y)F(y | x) dy \\
		& = F(\inner{R_x}{\theta_{0h,\lambda}} | x) + \int [F(\inner{R_x}{\theta_{0h,\lambda}} - hz | x) - F(\inner{R_x}{\theta_{0h,\lambda}} | x)] k(z)dz\\
		& = \tau + O(h^{s+1}).
	\end{align*}
	The last equation holds due to Theorem \ref{th-thlam} and Assumptions \ref{ass-f} and \ref{ass-K}. Let $K\hi *(z) = 2k(z) \bK(z) = \frac{d}{dz} \bK^2(z)$, thus $\int K\hi *(z) dz = \lim\limits_{z \rightarrow \infty} \bK^2(z) = 1$. For the same reason, we have
	\begin{align*}
		&\E[\bK^2_h(\inner{R_X}{\theta_{0h,\lambda}} - Y) | X = x]  \\
		& = h^{-1} \int K\hi* \left(\frac{\inner{R_X}{\theta_{0h,\lambda}} - y}{h}\right)F(y | x)dy \\
		& = \tau + O(h^{s+1}) + \int [F(\inner{R_x}{\theta_{0h,\lambda}} - hz | x) - F(\inner{R_x}{\theta_{0h,\lambda}} | x)] K\hi*(z)dz \\
		& = \tau + O(h^{s+1})  - h[f(\inner{R_x}{\theta_{0h,\lambda}} | x) + O(h)] \int z K\hi*(z)dz \\
		& = \tau + O(h^{s+1})  - h[f(\inner{R_x}{\theta_{0}} | x) + O(h) + O(h^{s+1})] \int z K\hi*(z)dz \\
		& = \tau - hB(x) \int z K\hi*(z)dz + O(h^2).
	\end{align*}
	By the proof of Theorem 3 of \cite{fernandes2019}, $\int z K\hi*(z)dz = 2 \int_0^{\infty} \bK(z) [1 - \bK(z)]dz < \infty$. In summary,
	$\E\{|\bK_h(-\eps(\theta_{0h,\lambda})) - \tau|^2 | X = x\}  = \tau - \tau^2 - 2hB(x) \int_0^{\infty} \bK(z) [1 - \bK(z)]dz + O(h^2)$. Since $\E[\norm{R_X}^2] = O(\omega^{-1})$ by Lemma S.4 of \cite{shang2015b}, 
	$E[\norm{S_{n,h\lambda}(\theta_{0h,\lambda})}^2]  \leq Cn^{-1}\tau(1-\tau) \E[\norm{R_X}^2] = O((n\omega)^{-1})$. This indicates that with probability approaching one, $\norm{S_{n,h\lambda}(\theta_{0h,\lambda})} \leq C_3(nw)^{-1/2}$ for some sufficiently large constant $C_3$. Let $r_{2n} = 4C_3(nw)^{-1/2}$ and $\mB (r_{2n}) = \{\theta \in \ca H: \norm{\theta} \leq r_{2n}\}$.

	Next we handle $I_1$. Let $T = (Y, X)$ denote the data variable. Define $A_n = \cap_{i=1}^n A_{ni} $, where
	$$
	A_{ni} = \{\|X_i\|_{L^2} \leq C\log n\},
	$$
	and $C$ is a positive constant satisfying $C\log n > 1$ and $C \geq \max(\sup_{u \in \Re} |k(u)|, \sup_{u \in \Re} |k^{(1)}(u)|)$. Because of Assumption \ref{ass-X}, we can choose a sufficiently large $C$ such that $P(A_n)$ converges to 1 as $n \rightarrow \infty$. Furthermore, we can even ensure $P(A_{ni}^c) = O(n^{-1})$ holds uniformly over $i = 1, \ldots, n$ for a enough large $C$. 
	
	Let $T = (Y, X)$ denote the data variable. 
	To deal with the first term $I_1$, we define a function
	$$
	\psi(T_i; \theta) = k_h(-\epsilon_i(\theta_{0h,\lambda})) \inner{R_{X_i}}{\theta} \bbone(A_{ni}). 
	$$
	For any $\theta_j = (\alpha_j, \beta_j) \in \ca H, j = 1, 2$, we have
	\begin{align*}
		|\psi(T_i; \theta_1) - \psi(T_i; \theta_2)| & =  |k_h(-\epsilon_i(\theta_{0h,\lambda})) | |\inner{R_{X_i}}{\theta_1 - \theta_2}| \bbone(A_{ni}) \\
		& \leq (Ch^{-1}) \left[|\alpha_1 - \alpha_2| + \|X_i\|_{L^2}\norm{\beta_1 - \beta_2}_{L^2}\right] \bbone(A_{ni}) \\
		& \leq (C^2h^{-1}\log n) \norm{\theta_1 - \theta_2}_2. 
	\end{align*}
	Let $\psi_n = (T_i; \theta) = C^{-2}(h^{-1}\log n)^{-1} \psi(T_i; \theta)$. Then $\psi_n$ satisfies that for any $\theta, \ti\theta \in \ca H$, 
	$$
	|\psi_n(T; \theta) - \psi_n(T; \ti\theta)| \leq \|{\theta - \ti\theta}\|_2. 
	$$
	
	By Lemma 3.1 of \cite{shang2015}, for any $\theta \in \ca H$, $\norm{\theta}_2 \leq \kappa \omega ^{-(2a+1)/2}\norm{\theta}$ for some constant $\kappa > 0$. Let $d_n = \ka \omega ^{-(2a+1)/2}$ and $p_n = d_n^{-2}\lambda^{-1}$. For any $\theta \in \ca H \setminus \{0\}$, let $\bar{\theta} = (\bar\alpha, \bar\beta) := \theta/(d_n\|\theta\|)$. Obviously $\|\bar\theta\|_2 \leq 1$. In other words, $|\bar\alpha| + \|\bar\beta\|_{L^2} \leq 1$. We can further show that $\lambda J(\bar\beta, \bar\beta) \leq \|\bar\beta\|^2 = d_n^{-2}$, which implies that $J(\bar\beta, \bar\beta) \leq d_n^{-2} \lambda^{-1} = p_n$. If we define $\ca F_{q_n} = \{\theta=(\alpha, \beta) \in \ca H: |\alpha| \leq 1, \|\beta\|_{L^2} \leq 1, J(\beta, \beta) \leq q_n\}$, where $q_n \geq 1$, 
	then $\bar\theta \in \ca F_{p_n}$. By Lemma 3.4 of \cite{shang2015}, we have for any $\theta \in \mB(r_{2n})$, with probability approaching one, 
	$$
	\norm{\sum_{i=1}^n [\psi_n(T_i; \bar\theta)R_{X_i} - \E\{\psi_n(T_i;\bar\theta)R_{X_i}\}]  } \leq C_5 (n^{1/2}p_n^{1/(4m)} + 1)(\omega^{-1}\log \log n)^{1/2},
	$$
	for some positive constant $C_5$. As a result,
	\begin{align*}
		&\norm{\sum_{i=1}^n [\psi_n(T_i; \bar\theta)R_{X_i} - \E\{\psi_n(T_i;\bar\theta)R_{X_i}\}]  }  \leq \\
		& ~~~~~C_5 \ka \omega^{-(a+1)}  (n^{1/2}p_n^{1/(4m)} + 1) (C^2h^{-1}\log n) (\log \log n)^{1/2} \|\theta\|
	\end{align*}
	with probability approaching 1. 
	On the other hand, by Cauchy-Schwarz inequality, we have
	\begin{align*}
		&\quad\quad \|\E[k_h(-\eps_i(\theta_{0h,\lambda}))\inner{R_{X_i}}{\theta}R_{X_i} \bbone(A^c_{ni})]\| \\
		& = \sup_{\|\theta_1\| = 1} | \E[k_h(-\eps_i(\theta_{0h,\lambda}))\inner{R_{X_i}}{\theta}\inner{R_{X_i}}{\theta_1} \bbone(A^c_{ni})]| \\
		& \leq Ch^{-1} \sup_{\|\theta_1\| = 1} \{\E[|\inner{R_{X_i}}{\theta}| \cdot |\inner{R_{X_i}}{\theta_1}|  \bbone(A^c_{ni})]\} \\
		& \leq Ch^{-1} \E[\inner{R_{X_i}}{\theta}^4]^{1/4} \sup_{\|\theta_1\| = 1}\E[\inner{R_{X_i}}{\theta_1}^4]^{1/4}  P(A^c_{ni})^{1/2}\\
		& \leq Ch^{-1} \|\theta\| P(A^c_{ni})^{1/2} = o(1) \|\theta\|, 
	\end{align*}
	where the last equation holds since $P(A_{ni}^c) = O(1/n)$ uniformly over $i  = 1, \ldots, n$ and $n^{-1/2}h^{-1} = o(1)$. 
	
	Now we are able to determine the order of $I_1$. Actually, with probability approaching 1, we have
	\begin{align*}
		&\quad\| DS_{n, h\lambda}(\theta_{0h,\lambda})\theta - DS_{h\lambda}(\theta_{0h,\lambda})\theta \| \\
		& \leq n^{-1} \norm{ \sum_{i = 1}^n [k_h(-\epsilon_i(\theta_{0h,\lambda})) \inner{R_{X_i}}{\theta}R_{X_i}\bbone(A_{ni}) - \E\{k_h(-\epsilon_i(\theta_{0h,\lambda})) \inner{R_{X_i}}{\theta} R_{X_i}\bbone(A_{ni})\}] } \\
		&\quad + n^{-1} \sum_{i = 1}^n \norm{\E[k_h(-\eps_i(\theta_{0h,\lambda}))\inner{R_{X_i}}{\theta}R_{X_i} \bbone(A^c_{ni})]} \\
		& = O_P\left(n^{-1/2}\omega^{-(a+1)-\frac{2l-2a-1}{4m}} (h^{-1} \log n) (\log \log n)^{1/2}\right)\|\theta\| + o(1)\|\theta\| \\
		& = o_P(1) \|\theta\|,
	\end{align*}
	where the last equation follows by the assumption \eqref{eq-lamorder}. Hence, with probability approaching one, we have $\| DS_{n, h\lambda}(\theta_{0h,\lambda})\theta - DS_{h\lambda}(\theta_{0h,\lambda})\theta \| \leq r_{2n}/18$ for any $\theta \in \mB(r_{2n})$. 
	
	Lastly, we treat $I_2$. Taylor's expansion leads to that
	\begin{align*}
		&\quad S_{n, h\lambda}(\theta_{0h, \lambda} + \theta) - S_{n, h\lambda}(\theta_{0h, \lambda}) - DS_{n, h\lambda}(\theta_{0h, \lambda})\theta \\
		& = \int_0^1 \int_0^1 D^2S_{n,h\lambda}(\theta_{0h,\lambda} + s^{\prime} s\theta) (s\theta) \theta ds^{\prime}ds. 
	\end{align*}
	We assume that $\cap_{i = 1}^n A_{ni}$ holds for the rest of the proof. For an arbitrary $\theta \in \ca H \setminus \{0\}$, let $\bar\theta = \theta/(d_n\|\theta\|)$, so
	$\bar\theta \in \ca F \lo {p_n}$. Let
	$$
	\psi_n(T_i; \theta) = \frac{\|R_{X_i}\| \cdot \inner{R_{X_i}}{\bar\theta}}{\sqrt{2C_R} (C\log n)^2 \omega^{-(2a+1)/2}} \bbone(A_{ni}), i = 1, \ldots, n,
	$$
	where $C_R$ is a constant such that $\|R_{X_i}\|^2 \leq C_R[1 + (C\log n)^2 \omega^{-(2a+1)}]$ based on Lemma S.4 of \cite{shang2015b}. Hence $\psi_n$ satisfies Lipschitz continuity given in (3.3) of \cite{shang2015}. Then by Lemma 3.4 of \cite{shang2015}, we have, with probability approaching 1, for any $\theta \in \ca H \setminus \{0\}$,
	$$
	n^{-1/2} \norm{ \sum_{i = 1}^n [\psi_n(T_i; \bar\theta) R_{X_i} - \E_T\{\psi_n(T_i; \bar\theta)R_{X_i}\} ]} \leq C_6 p_n^{1/(4m)} (\omega^{-1}\log \log n)^{1/2}, 
	$$
	for some positive constant $C_6$. This implies that with probability approaching one, 
	\begin{align*}
		& \quad \norm{ \sum_{i = 1}^n [\| R_{X_i}\| \cdot  \inner{R_{X_i}}{\theta} \bbone(A_{ni}) R_{X_i}  - \E_T \{\| R_{X_i}\| \cdot  \inner{R_{X_i}}{\theta} \bbone(A_{ni}) R_{X_i} \} ] } \\
		& \leq C_7 n^{1/2} \omega^{-(2a+\frac{3}{2})-\frac{2l-2a-1}{4m}} (\log n)^2 (\log \log n)^{1/2} \|\theta\|
	\end{align*}
	for some large constant $C_7$. Note that this inequality also holds for $\theta = 0$. On the other hand, by Lemma S.3 and S.4 of \cite{shang2015b} and Cauchy inequality,
	$$
	\E_T\{\|R_{X}\| \cdot |\inner{R \lo {X}}{\theta}|^2\} \leq (\E_T \{\{\|R_{X}\|^2\})^{1/2} \cdot (\E_T \{\{|\inner{R \lo {X}}{\theta}|^4\})^{1/2} \leq C_8\omega^{-1/2} \|\theta\|^2
	$$
	for some large constant $C_8$. Therefore, on $\cap_{i = 1}^nA_{ni}$, for any $\theta \in \mB(r_{2n})$ and some $s, s^{\prime} \in [0, 1]$, we have
	\begin{align*}
		& \quad \|D^2S_{n,h\lambda}(\theta_{0h,\lambda} + s^{\prime} s\theta) \theta \theta ds^{\prime} \| \\
		& = n^{-1} \norm{ \sum_{i = 1}^n h^{-2}k^{(1)}(-\eps_i(\ti\theta)/h)\inner{R \lo {X_i}}{\theta}^2  {R \lo {X_i}} } \\
		& \leq n^{-1} h^{-2} \sup_{a \in \Re} |k^{(1)}(u)| \sum_{i = 1}^n \|R_{X_i}\| \cdot |\inner{R \lo {X_i}}{\theta}|^2 \\
		& \leq Cn^{-1} h^{-2} \inner{\sum_{i = 1}^n  [\|R_{X_i}\| \cdot \langle R_{X_i}, \theta \rangle \bbone(A_{ni}) R_{X_i} - \E_T\{\|R_{X_i}\| \cdot \inner{R \lo {X_i}}{\theta}  \bbone(A_{ni}) R_{X_i}\}]} {\theta} \\
		& \quad + Ch^{-2} \E_T \{\|R_{X_i}\| \cdot |\inner{R \lo {X_i}}{\theta}|^2  \bbone(A_{ni})\} \\
		& \leq C^{\prime} [n^{-1}\omega^{-2(a+1)-\frac{2l-2a-1}{4m}} (\log n)^2 h^{-2} (\log \log n)^{1/2} + n^{-1/2}h^{-2} \omega^{-1} ] \times \|\theta\| \\
		& \leq \|\theta\|/18,
	\end{align*}
	where the last inequality is obtained from assumption \eqref{eq-lamorder}. Combining the above results, we have that, with probability approaching one, for any $\theta \in \mB(r_{2n})$, 
	$$
	T_2(\theta) \leq \|I_1\| + \|I_2\| + \|I_3\| \leq 2(r_{2n}/4 + r_{2n}/18 + r_{2n}/18) = 13r_{2n}/18. 
	$$
	That is to say, $T_2(\mB(r_{2n})) \subset \mB(r_{2n})$ with probability approaching to 1. 
	
	Next we need to show that $T_2$ is a contraction mapping. Given any $\theta_1, \theta_2 \in \mB(r_{2n})$, we have
	\begin{align*}
		&\quad T_2(\theta_1) - T_2(\theta_2) \\
		& = \theta_1 - \theta_2 - [DS_{h\lambda}(\theta_{0h,\lambda})]^{-1} [S_{n,h\lambda}(\theta_{0h,\lambda} + \theta_1) - S_{n,h\lambda}(\theta_{0h,\lambda} + \theta_2)] \\
		& = - [DS_{h\lambda}(\theta_{0h,\lambda})]^{-1} \int_0^1 \int_0^1 D^2S_{n,h\lambda}(\theta_{0h,\lambda} + s^{\prime}(\theta_2 + s(\theta_1 - \theta_2))) (\theta_2 + s(\theta_1 - \theta_2))(\theta_1 - \theta_2)ds^{\prime}ds \\
		& \quad - [DS_{h\lambda}(\theta_{0h,\lambda})]^{-1} [DS_{n,h\lambda}(\theta_{0h,\lambda}) - DS_{h\lambda}(\theta_{0h,\lambda})](\theta_1 - \theta_2) \\
		& := I_4 + I_5. 
	\end{align*}
	Using the same arguments as in the analysis of the terms $I_2$ and $I_1$, we can show that with probability approaching one, 
	\begin{align*}
		\|I_4\| & \leq O(n^{-1}\omega^{-2(a+1)-\frac{2l-2a-1}{4m}} (\log n)^2 h^{-2} (\log \log n)^{1/2} + n^{-1/2}h^{-2} \omega^{-1}) \|\theta_1 - \theta_2\| \\
		& \leq \|\theta_1 - \theta_2\| /3,
	\end{align*}
	and
	\begin{align*}
		\|I_5\| & = O\left(n^{-1/2}\omega^{-(a+1)-\frac{2l-2a-1}{4m}} (h^{-1} \log n) (\log \log n)^{1/2}\right) \|\theta_1 - \theta_2\| + o(1) \|\theta_1 - \theta_2\| \\
		& \leq \|\theta_1 - \theta_2\|/3.
	\end{align*}
	It follows that $\|T_2(\theta_1) - T_2(\theta_2) \| \leq 2\|\theta_1 - \theta_2\|/3$ with probability approaching 1. Therefore, $T_2$ is a contraction mapping from $\mB(r_{2n})$ to itself.
	By the contraction mapping theorem, there must exist a unique element $\ti\theta \in \mB(r_{2n})$ such that $T_2(\ti\theta) = \ti\theta$, which implies that $S_{n,h\lambda}(\theta_{0h,\lambda} + \ti\theta) = 0$. Let $\hat{\theta}_{h\lambda} = \theta_{0h,\lambda} + \ti\theta$, then $S_{n,h\lambda}(\hat{\theta}_{h\lambda}) = 0$. That is to say, $\hat{\theta}_{h\lambda}$ is the minimizer of the loss function $\ell_{n,h\lambda}(\theta)$. Furthermore, by Theorem \ref{th-thlam}, with probability approaching one,
	$$
	\|\hat{\theta}_{h\lambda} - \theta_0\| \leq \|\hat{\theta}_{h\lambda}  -  \theta_{0h,\lambda}\| + \| \theta_{0h,\lambda} - \theta_0\| \leq r_{2n} + r_{1n} = O((n\omega)^{-1/2} + \omega^l)= O(r_n). 
	$$
	This completes the proof.

	\subsection{Proof of Theorem \ref{th-Bahadur}}

	By Theorem \ref{th-rate}, there exists some sufficiently large $M > 0$ such that, with probability approaching one, $\|\hat{\theta}_{h,\lambda} - \theta_0 \| \leq Mr_n$. Let $\theta = \hat{\theta}_{h,\lambda} - \theta_0$ for notational convenience. In the rest of the proof we assume that $\|\theta\| \leq Mr_n$ because the probability of its complement is negligible. Let $d_n = \ka M \omega^{-(2a+1)/2}r_n$ and $\ti\theta = d_n^{-1}\theta$. Let $p_n = \ka^{-2} \omega^{1-2m}$, where $\ka$ is defined in the proof of Theorem \ref{th-rate}. Note that $p_n \geq 1$ for sufficiently large $n$ as $\omega$ decays to 0 and $1 - 2m = 1 - 2(k - a) < 0$. We can show that $\|\theta\| \leq M r_n$ indicates that $\ti\theta \in \ca F _{p_n}$. Actually, by Lemma 3.1 of \cite{shang2015}, $\|\ti\theta\|_2 = d_n^{-1}\|\theta\|_2 \leq d_n^{-1}\ka \omega^{-(2a+1)/2} \|\theta\| \leq d_n^{-1}\|\theta\|_2 \leq d_n^{-1}\ka \omega^{-(2a+1)/2}Mr_n = 1$. Thus $|\ti\alpha| \leq 1$ and $\|\ti\beta\|_{L^2} \leq 1$. Additionally, 
	$$
	J(\ti\beta, \ti\beta) = d_n^{-2}\lambda^{-1}(\lambda J(\beta,\beta)) \leq d_n^{-2}\lambda^{-1} \|\theta\|^2 \leq d_n^{-2}\lambda^{-1}(Mr_n)^2 = \ka^{-2}\omega^{1-2m} = p_n. 
	$$
	Therefore, $\ti\theta \in \ca F _{p_n}$. 
	
	For $i = 1, \ldots,n $, define
	$
	A_{ni} = \{\|X_i\|\lo{L^2} \leq C\log n\},
	$
	where $C$ is a positive constant satisfying $C\log n > 1$ and $C \geq \sup_{u \in \Re} |k(u)|$. Denote $\cap_{i=1}^n A_{ni} $ by $A_n$. 
	Because of Assumption \ref{ass-X}, we can choose a sufficiently large $C$ such that $P(A_n)$ converges to 1 as $n \rightarrow \infty$. Furthermore, we can even ensure that 
	$$
	n^{1/2}\omega^{-1/2}P(\bbone(A_{ni}^c))^{1/4} = o(p_n^{1/(4m)}(\omega^{-1}\log\log n)^{1/2})
	$$ 
	for a enough large $C$. 
	
	Let $D_n = (C^2\log n \cdot h^{-1})^{-1}d_n^{-1}$. We define $\psi(T_i; \theta) = \bK(-\eps_i(\theta + \theta_0)/h) - \bK(-\eps_i(\theta_0)/h)$ and $\psi_n(T_i; \ti\theta) = D_n\psi(T_i; d_n\ti\theta) \bbone(A_{ni})$. Then for any $\ti\theta_1, \ti\theta_2 \in \ca F_{p_n}$, we have
	\begin{align*}
		|\psi_n(T_i; \ti\theta_1) - \psi_n(T_i; \ti\theta_2) & = D_n|\psi(T_i; d_n\ti\theta_1) - \psi(T_i; d_n\ti\theta_2)| \bbone(A_{ni}) \\
		& = D_n |\bK(-\eps_i(d_n\ti\theta_1 + \theta_0)/h) - \bK(-\eps_i(d_n\ti\theta_2 + \theta_0)/h)|\bbone(A_{ni}) \\
		& \leq D_n \sup_{u \in \Re}|k(u)| \cdot h^{-1} \cdot d_n \cdot |\inner{R_{X_i}}{\ti\theta_1 - \ti\theta_2}| \bbone(A_{ni}) \\
		& \leq D_n (Ch^{-1}) d_n (C\log n) \|\ti\theta_1 - \ti\theta_2\|_2 =  \|\ti\theta_1 - \ti\theta_2\|_2. 
	\end{align*}
	It follows by Lemma 3.4 of \cite{shang2015}, that there exists a constant $C_3$ such that for sufficiently large $n$, with probability approaching one, 
	\begin{align*}
		\|H_n(\ti\theta)\| & = \frac{1}{\sqrt{n}} \norm{\sum_{i = 1}^n [\psi_n(T_i; \ti\theta)R_{X_i} - \E_T\{\psi_n(T; \ti\theta)R_{X}\}] } \\
		& \leq C_3(p_n^{1/(4m)}\|\ti\theta\|_2^{\gamma} + n^{-1/2})(\omega^{-1}\log \log n)^{1/2} \\
		& \leq C_3(p_n^{1/(4m)} + n^{-1/2})(\omega^{-1}\log \log n)^{1/2},
	\end{align*}
	where $\gamma = 1 - 1/(2m)$. On the other hand, by Cauchy's inequality and the mean value theorem, we have
	\begin{align*}
		& \quad	\|\E_T\{\psi(T_i; d_n\ti\theta)R_{X_i} \bbone(A_{ni}^c) \} \| \\
		& \leq \E_T \{|\psi(T_i; d_n\ti\theta)| \cdot \|R_{X_i}\| \cdot \bbone(A_{ni}^c) \} \\
		& \leq \E_T \{\sup_{s \in [0,1]} |k_h(-\eps(\theta_0 + sd_n\ti\theta))| d_n |\inner{R_{X_i}}{\ti\theta}| \cdot \|R_{X_i}\| \cdot \bbone(A_{ni}^c)\} \\
		& \leq C_4 d_n \E\{(1 + \|X_i\|_{L^2}) \|R_{X_i}\| \bbone(A_{ni}^c)\} \\
		& \leq C_4 d_n \E_T\{(1 + \|X_i\|_{L^2})^4\}^{1/4} P(A_{ni}^c)^{1/4} \E\{ \|R_{X_i}\| ^2\}^{1/2} 
	\end{align*}
	As a result, based on the choice of $C$, there exists some constant $C^{\prime} > 0$, 
	\begin{align*}
		&\quad D_n n^{1/2} \|\E_T\{\psi(T_i; d_n\ti\theta)R_{X_i} \bbone(A_{ni}^c) \} \| \\
		& \leq C^{\prime} \E_T\{(1 + \|X_i\|_{L^2})^4\}^{1/4} n^{1/2} \omega^{-1/2}  P(A_{ni}^c)^{1/4}  (\log n)^{-1} h \\
		& \leq C^{\prime} p_n^{1/(4m)} (\omega^{-1} \log \log n)^{1/2}
	\end{align*}
	Note that $\bbone(A_{ni}) = 1$ on $A_n$. Therefore, with probability approaching one, 
	\begin{align*}
		& \quad n^{-1/2}D_n \norm{ \sum_{i = 1}^n [\psi(T_i; \theta)R_{X_i} - \E_T\{\psi(T; \theta)R_{X}\}]} \\
		& = \| n^{-1/2} \sum_{i = 1}^n [\psi_n(T_i; \ti\theta)R_{X_i} - \E_T\{\psi_n(T; \ti\theta)R_{X}\}] - D_n n^{1/2} \E_T\{\psi(T_i; d_n\ti\theta)R_{X_i} \bbone(A_{ni}^c) \} \|\\
		& \leq \|H_n(\ti\theta)\| + D_n n^{1/2} \|\E_T\{\psi(T_i; d_n\ti\theta)R_{X_i} \bbone(A_{ni}^c) \} \| \\
		& \leq C_5 p_n^{1/(4m)} (\omega^{-1} \log \log n)^{1/2}. 
	\end{align*}
	Then with probability approaching one, we have
	\begin{align*}
		& \quad \|S_{n,h\lambda}(\theta + \theta_0) - S_{n,h\lambda}(\theta_0) - \E_T\{S_{n,h\lambda}(\theta + \theta_0) - S_{n, h\lambda}(\theta_0)\} \| \\
		& = n^{-1} \| \sum_{i = 1}^n [\psi(T_i; \theta)R_{X_i} - \E_T\{\psi(T_i; \theta)R_{X_i}\}] \| \\
		& \leq C_5 n^{-1/2} p_n^{1/(4m)} D_n^{-1} (\omega^{-1} \log \log n)^{1/2} \\
		& = C_5 \ka^{\gamma} C^2 M n^{-1/2} \omega^{-\frac{4ma+6m-1}{4m}} r_n( h^{-1} \log n) (\log \log n)^{1/2}.
	\end{align*}
	
	Note that $\theta$ satisfies $S_{n,h\lambda}(\theta + \theta_0) = 0$. Additionally, from the proof of Theorem \ref{th-thlam} we know that $\int \rho_{\tau}(u)D_{\theta}^2 f_{\eps}(u; \theta_0) \theta du + P_{\lambda}(\theta) = \theta$. Recall that $\ti{K}_h(u; \theta) = \E[k_h(u - \eps(\theta))]$, where the expectation is taken with respect to the distribution of $\eps(\theta)$. 
	It follows that
	\begin{align*}
		& \quad \E_T\{S_{n,h\lambda}(\theta + \theta_0) - S_{n, h\lambda}(\theta_0)\} \\
		& = S_{h\lambda}(\theta + \theta_0) - S_{h\lambda}(\theta_0) \\
		& = P_{\lambda}(\theta) + \int \rho_{\tau}(u)  D_{\theta} \ti{K}_h(u; \theta + \theta_0) du - \int \rho_{\tau}(u)  D_{\theta} \ti{K}_h(u; \theta_0) du \\
		& = P_{\lambda}(\theta) + \int \rho_{\tau}(u)D_{\theta}^2 f_{\eps}(u; \theta_0) \theta du - \int \rho_{\tau}(u)D_{\theta}^2 f_{\eps}(u; \theta_0) \theta du \\
		& \quad + \int \rho_{\tau}(u)  D_{\theta} \ti{K}_h(u; \theta + \theta_0) du - \int \rho_{\tau}(u)  D_{\theta} \ti{K}_h(u; \theta_0) du \\
		& = \theta + \int \rho_{\tau}(u)  D_{\theta} \ti{K}_h(u; \theta + \theta_0) du - \int \rho_{\tau}(u)  D_{\theta} \ti{K}_h(u; \theta_0) du - \int \rho_{\tau}(u)D_{\theta}^2 f_{\eps}(u; \theta_0) \theta du \\
		& = \theta + \int \rho_{\tau}(u) \left[ \int_0^1 D^2_{\theta}\ti{K}_h(u;  \theta_0 + s\theta)\theta ds\right] du - \int \rho_{\tau}(u)D_{\theta}^2 f_{\eps}(u; \theta_0) \theta du \\
		& = \theta  -  \int \rho_{\tau}(u) \left(\int_0^1[D^2_{\theta}f_{\eps}(u; \theta_0)\theta - D_{\theta}^2 \ti{K}_h(u; s\theta + \theta_0) \theta]ds\right)du 
	\end{align*}
	From the proof of Theorem \ref{th-thlam} and $h = O(\omega^{k/(s+1)})$, we know that  
	$$\norm{ \int \rho_{\tau}(u) \left(\int_0^1[D^2_{\theta}f_{\eps}(u; \theta_0)\theta - D_{\theta}^2 \ti{K}_h(u; s\theta + \theta_0) \theta]ds\right)du }      \leq C[o(h^s) + \|\theta\|] \cdot \|\theta\| \leq 2h^{-1} r_n^2$$ 
	uniformly over $\tau \in [\underline{\tau}, \bar{\tau}]$ for sufficiently large $n$.  Therefore, as $n \rightarrow \infty$, with probability approaching one, 
	\begin{align*}
		\|\theta + S_{n, h\lambda}(\theta_0) \| \leq C_5 \ka^{\gamma} C^2 M n^{-1/2} \omega^{-\frac{4ma+6m-1}{4m}} r_n( h^{-1} \log n) (\log \log n)^{1/2} + 2h^{-1}r_n^2. 
	\end{align*}
	This completes the proof. 
	
	\subsection{Proof of Corollary \ref{cor-pointwise}}

	Let $\zeta_i = \bK({-\eps_i(\theta_0)}/{h}) - \tau, i = 1, \ldots, n$. 
	By the definition of $K_t$, we have 
	\begin{align*}
		\|K_t\|_1^2 & = \inner{\sum_{\nu} \frac{\phi_{\nu}(t)}{1 + \lambda \rho_{\nu}} \phi\lo{\nu}}  {\sum_{\nu} \frac{\phi_{\nu}(t)}{1 + \lambda \rho_{\nu}}\phi\lo{\nu}}_1 \\
		& 	\leq C_{\phi}^2 \sum_{\nu} \frac{\nu^{2a}}{1 + \lambda\rho_{\nu}} \leq 	\sum_{\nu}\frac{{\nu}^{2a}}{1 + c_{\rho}\lambda {\nu}^{2l}}  \\
		& \leq C_K \omega^{-(2a+1)},
	\end{align*}
	where $C_K$ is a constant that only depends on $C_{\phi}$ and $c_{\rho}$.
	Then it follows that
	\begin{align*}
		&\quad |\hat{\beta}_{h,\lambda}(t) - \beta_0(t, \tau) + \frac{1}{n} \sum_{i = 1}^n \zeta_i \eta(X_i)(t) + (W_{\lambda}\beta_0(\cdot, \tau))(t)| \\
		& = \inner{K_t}{\hat{\beta}_{h,\lambda} - \beta_0(\cdot, \tau) +  \frac{1}{n} \sum_{i = 1}^n \zeta_i \eta(X_i) + (W_{\lambda}\beta_0(\cdot, \tau))} _1 \\
		& \leq \| K_t \|_1 \norm{\hat{\beta}_{h,\lambda} - \beta_0(\cdot, \tau) +  \frac{1}{n} \sum_{i = 1}^n \zeta_i \eta(X_i) + (W_{\lambda}\beta_0(\cdot, \tau))} _1 \\
		& = O_P(a_n\omega^{-(2a+1)/2}),
	\end{align*}
	where the last equation is obtained from Theorem \ref{th-Bahadur}. 
	
	Let $\bK_h(\cdot)= \bK(\cdot/h)$. We first employ the law of iterated expectations to find $\E \zeta_i^2 [\eta(X_i)(t)]^2$. In particular, given $X = x$, 
	\begin{align*}
		\E\{|\bK_h(-\eps(\theta_{0})) - \tau|^2 | X = x\} & = \E[\bK^2_h(\inner{R_X}{\theta_{0}} - Y) | X = x]  \\
		&\quad -2\tau \E[\bK_h(\inner{R_X}{\theta_{0}} - Y) | X = x]  + \tau^2.
	\end{align*}
	Integration by parts leads to 
	\begin{align*}
		& \E[\bK_h(\inner{R_X}{\theta_{0}} - Y) | X = x] \\
		& = \int \bK_h(\inner{R_x}{\theta_{0}} - y) f(y | x) dy \\
		& = \int k_h(\inner{R_x}{\theta_{0}} - y)F(y | x) dy \\
		& = F(\inner{R_x}{\theta_{0}} | x) + \int [F(\inner{R_x}{\theta_{0}} - hz | x) - F(\inner{R_x}{\theta_{0}} | x)] k(z)dz\\
		& = \tau + O(h^{s+1}).
	\end{align*}
	The last equation holds due to the definition of $\theta(\tau)$ and Assumptions \ref{ass-f} and \ref{ass-K}. Let $K\hi *(z) = 2k(z) \bK(z) = \frac{d}{dz} \bK^2(z)$, thus $\int K\hi *(z) dz = \lim\limits_{z \rightarrow \infty} \bK^2(z) = 1$. For the same reason, we have
	\begin{align*}
		&\E[\bK^2_h(\inner{R_X}{\theta_{0}} - Y) | X = x]  \\
		& = h^{-1} \int K\hi* \left(\frac{\inner{R_X}{\theta_{0}} - y}{h}\right)F(y | x)dy \\
		& = \tau + O(h^{s+1}) + \int [F(\inner{R_x}{\theta_{0}} - hz | x) - F(\inner{R_x}{\theta_{0}} | x)] K\hi*(z)dz \\
		& = \tau + O(h^{s+1})  - h[f(\inner{R_x}{\theta_{0}} | x) + O(h)] \int z K\hi*(z)dz \\
		& = \tau + O(h^{s+1})  - h[f(\inner{R_x}{\theta_{0}} | x) + O(h) + O(h^{s+1})] \int z K\hi*(z)dz \\
		& = \tau - hB(x) \int z K\hi*(z)dz + O(h^2).
	\end{align*}
	By the proof of Theorem 3 of \cite{fernandes2019}, $\int z K\hi*(z)dz = 2 \int_0^{\infty} \bK(z) [1 - \bK(z)]dz < \infty$. In summary,
	$\E\{|\bK_h(-\eps(\theta_{0})) - \tau|^2 | X = x\}  = \tau - \tau^2 - 2hB(x) \int_0^{\infty} \bK(z) [1 - \bK(z)]dz + O(h^2)$ uniformly for $x$. 
	Therefore, by Assumption \ref{ass-f} (c), 
	\begin{align*}
		\E \zeta_i^2 [\eta(X_i)(t)]^2 & \leq C_1 \E \zeta_i^2 B(X_i) [\eta(X_i)(t)]^2 \\
		& \asymp C_1 (\tau - \tau^2) \sum_{\nu} \frac{\E(B(X_i) X_{i\nu}^2) }{(1 + \lambda\rho_{\nu})^2} \phi_{\nu}^2(t) \\
		& \asymp \sum_{\nu} \frac {\phi_{\nu}^2(t)} {(1 + \lambda\rho_{\nu})^2} \geq \sigma_t^2 \omega^{-(2a + 1)}
	\end{align*}
	for some positive constant $\sigma_t^2$. 
	Next we find the $\E \zeta_i [\eta(X_i)(t)]$. It is actually $R^{(1)}_h(\theta_0; \tau)\theta_1$ if we take $\theta_1 = (0, K_t)$. 
	By the proof of (ii) of Lemma \ref{lemma1}, we have
	\begin{align*}
		|R^{(1)}_h(\theta_0; \tau)\theta_1 - R^{(1)}(\theta_0; \tau)\theta_1| & \leq Ch^{s+1} \E[|\inner{R \lo {X}}{ \theta_1}|]  \\
		& \leq  Ch^{s+1} \left[\E\left(\langle R_X, \theta_1 \rangle^4\right)\right]^{1/4} \\
		& \leq C h^{s+1} \|\theta_1\| = C h^{s+1} \|K_t\|_1
	\end{align*}
	by the Cauchy-Schwarz inequality and Lemma S.3 of \cite{shang2015b}. On the other hand, as we've shown above, $\|K_t\|_1 \leq C_K \omega^{-(2a+1)/2}$.
	Since $R^{(1)}(\theta_0; \tau)\theta_1 = 0$, $\{\E \zeta_i [\eta(X_i)(t)]\}^2 = O(h^{2s+2}) \omega^{-(2a + 1)}$. 
	On the other hand, since $h = o(1)$, $\{\E\zeta_i [\eta(X_i)(t)]\}^2= o(\E \zeta_i^2 [\eta(X_i)(t)]^2)$. Hence we conclude that $s_n^2 \asymp n\omega^{-(2a + 1)}$. 
	
	We check the Lindeberg's condition to establish the central limit theorem for the triangular array $\sum_{i = 1}^n \{\zeta_i [\eta(X_i)(t)] - R^{(1)}_h(\theta_0; \tau)\theta_1 \}$, where $\theta_1 = (0, K_t) \in \ca H$. Rearranging $\zeta_i [\eta(X_i)(t)] - R^{(1)}_h(\theta_0; \tau)\theta_1$ leads to
	\begin{align*}
		\zeta_i [\eta(X_i)(t)] - R^{(1)}_h(\theta_0; \tau)\theta_1 & = \left[\int_{-\infty}^{\inner{R \lo {X_i}}{ \theta_0}} k_h(u - Y_i) du - \tau\right] \inner{R \lo {X_i}}{ \theta_1} - \\
		& \quad\quad \int \left\{\int_{-\infty}^{\inner{R \lo {x}}{ \theta_0}} (\E[k_h(u - Y | x)]) du - \tau\right\} \inner{R \lo {x}}{ \theta_1} d P_X(x) .
	\end{align*}
	Obviously, there exists a positive constant $C$ that is independent of $n$ such that $|{\zeta}_i| \leq C$ a.s. As shown above, $|\inner{R \lo {X}}{ \theta_1}| \leq \|\eta(X)\|_1 \|K_t\|_1 \leq C_K\omega^{-(2a+1)/2} \|\eta(X)\|_1$.
	Furthermore, by Lemma S.4 of \cite{shang2015b}, $\|\eta(X)\|_1^2 \leq C_R\|X\|_{L^2}^2 \omega^{-(2a+1)}$ for some constant $C_R > 0$. 
	
	For any $\epsilon > 0$, we choose $\ti{C}$ sufficiently large such that $\ti{C}\eps \sigma_t / C_K \log(1/\omega) > C$. Considering the fact that $[R^{(1)}_h(\theta_0; \tau)\theta_1]^2 = o(s_n^2/n)$, we can 
	obtain that
	\begin{align*}
		& \quad	ns_n^{-2} \E \{[\zeta \inner{R \lo {X}}{ \theta_1} - R^{(1)}_h(\theta_0; \tau)\theta_1]^2 \bbone([\zeta_i \inner{R \lo {X}}{ \theta_1} - R^{(1)}_h(\theta_0; \tau)\theta_1]^2 \geq \eps^2s_n^2)  \} \\
		& \lesssim ns_n^{-2} \E    \{\zeta^2 \inner{R \lo {X}}{ \theta_1}^2  \bbone([\zeta \inner{R \lo {X}}{ \theta_1}]^2 \geq \eps^2s_n^2)  \}  + ns_n^{-2} \E \{[R^{(1)}_h(\theta_0; \tau)\theta_1]^2  \bbone([\zeta\inner{R \lo {X}}{ \theta_1}]^2 \geq \eps^2s_n^2)  \} \\
		& \leq ns_n^{-2} C^2 \E\{\inner{R \lo {X}}{ \theta_1}^2 \bbone(\zeta^2\inner{R \lo {X}}{ \theta_1}^2 \geq \eps^2s_n^2)  \} + ns_n^{-2} \E \{[R^{(1)}_h(\theta_0; \tau)\theta_1]^2  \bbone([\zeta\inner{R \lo {X}}{ \theta_1}]^2 \geq \eps^2s_n^2)  \} \\
		& \leq O(\omega^{-(2a+1)}) (\E\{\inner{R \lo {X}}{ \theta_1}^4\})^{1/2} \cdot P(|\zeta| \cdot |\inner{R \lo {X}}{ \theta_1}| \geq \eps s_n)^{1/2} \\
		& \leq O(\omega^{-2(2a+1)}) P(|\zeta| \cdot \|\eta(X)\|_1 C_K\omega^{-(2a+1)/2} \geq \eps \sqrt{n} \sigma_t \omega^{-(2a + 1)/2})^{1/2} \\
		& =  O(\omega^{-2(2a+1)}) P(|\zeta| \cdot \|\eta(X)\|_1 \geq \eps \sqrt{n} \sigma_t/C_K)^{1/2} \\
		& \leq O(\omega^{-2(2a+1)}) \left(\left[P(|\zeta| \geq (\eps \sigma_t / C_K) \ti{C} \log(\omega^{-1}) )\right] ^{1/2} +  \left[P(\|\eta(X)\|_1 \geq \sqrt{n} /(\ti{C} \log(\omega^{-1}))\right]^{1/2}  \right) \\
		& \leq O(\omega^{-2(2a+1)}) \E{\exp(s\|X\|_{L^2})} \exp \left( -\frac{s\sqrt{n\omega^{2a+1}}} {C_R^{1/2} \ti{C} \log (\omega^{-1}) } \right) \\
		& = o(1),
	\end{align*}
	where the last equation holds because $n\omega^{2a+1}(\log(\omega^{-1}))^{-4} \rightarrow \infty$ as $n \rightarrow \infty$. Then by the Lindberg's CLT, we have
	$s_n^{-1} \sum_{i = 1}^n \{\zeta_i \eta(X_i)(t) - \E[\zeta_i \eta(X_i)(t)]\} \convd N(0, 1)$. Moreover, $s_n^{-1}n   \E[\zeta_i \eta(X_i)(t)] = O(\sqrt{n} \omega^{(2a + 1)/2} \cdot h^{s + 1} \omega^{-(2a + 1)/2}) = O(\sqrt{n} h^{s + 1}) = o(1)$. Therefore, following the condition that $n a_n \omega^{-(2a + 1)/2}/s_n \asymp n^{1/2}a_n = o(1)$, we get that as $n \rightarrow \infty$, 
	\begin{align*}
		& \quad \frac{n}{s_n}[\hat{\beta}_{h,\lambda} - \beta_0(t, \tau) + (W_{\lambda}\beta_0(\cdot, \tau))(t)]  \\
		& = s_n^{-1} \sum_{i = 1}^n \{-\zeta_i \eta(X_i)(t) + \E[\zeta_i \eta(X_i)(t)]\} - ns_n^{-1}  \E[\zeta_i \eta(X_i)(t)]  + O_P(n a_n \omega^{-(2a + 1)/2}/s_n ) \\
		& = -s_n^{-1} \sum_{i = 1}^n \{\zeta_i \eta(X_i)(t) - \E[\zeta_i \eta(X_i)(t)]\}  + o(1) + o_P(1) \convd N(0, 1). 
	\end{align*}
	
	\subsection{Proof of Theorem \ref{th-weak}}

	By the proof of Corollary \ref{cor-pointwise}, we have that for an arbitrary $\tau \in \mI$, 
	$$
	n^{1/2}\omega^{a + 1/2} \sup_{t \in \mI} |\hat{\beta}_{h,\lambda}(t, \tau) - \beta_0(t, \tau) + \ti{S}_n(\theta_0)(t) + (W_{\lambda}\beta_0(\cdot, \tau))(t)| = o_P(1),
	$$
	where 
	$$
	\ti{S}_n(\theta_0)(t) = \frac{1}{n} \sum_{i = 1}^n \left[\int_{-\infty}^{\inner{R \lo {X_i}}{ \theta_0}} k_h(u - Y_i) du - \tau\right] \inner{R \lo {X_i}}{ \theta_1}  - R^{(1)}_h(\theta_0; \tau)\theta_1
	$$
	with $\theta_1 = (0, K_t)$. Since $\beta_0 = \sum_{\nu} b_{\nu}\phi_{\nu}$, by Cauchy-Schwarz inequality, we have
	\begin{align*}
		(W_{\lambda}\beta_0(\cdot, \tau))(t) & = \sum_{\nu} b_{\nu} \frac{\lambda\rho_{\nu}}{1 + \lambda\rho_{\nu}} \phi_{\nu}(t) \\
		& \leq \lambda \left(\sum_{\nu} b_{\nu}^2\rho_{\nu}^2 \right)^{1/2} \cdot \left(\sum_{\nu} \frac {\phi_{\nu}^2(t)} {(1 + \lambda\rho_{\nu})^2} \right)^{1/2} \\
		& \asymp \lambda \omega^{-(a + 1/2)}
	\end{align*}
	uniformly over $t \in \mI$. It follows by the condition $n^{1/2}\lambda = o(1)$ that
	\begin{equation} \label{eq-strongapp}
		n^{1/2}\omega^{a + 1/2} \sup_{t \in \mI} |\hat{\beta}_{h,\lambda}(t, \tau) - \beta_0(t, \tau) + \ti{S}_n(\theta_0)(t)| = o_P(1),
	\end{equation}

	To show weak convergence of $n^{1/2}\omega^{a + 1/2}[\hat{\beta}_{h,\lambda}(\cdot, \tau) - \beta_0(\cdot, \tau)]$ to some Gaussian process in $H^m(\mI)$, we only need to verify that $J_n(t) := n^{1/2}\omega^{a + 1/2} \ti{S}_n(\theta_0)(t)$ converges weakly to the Gaussian process $\ca G$ in $H^m(\mI)$ equipped with the inner product $V(\cdot, \cdot)$. Note that $\{\phi_{\nu}\}$ is a complete orthonormal basis of $H^m(\mI)$ by Assumption \ref{ass-diag}. Then by Theorem 1.8.4 of \cite{van1996weak}, it is equivalent to prove that (i) $J_n(\cdot)$ is asymptotically finite-dimensional and (ii) $V(J_n, \phi_{\nu})$ converges in distribution to $V(\ca G, \phi_{\nu})$ for each $\nu \geq 1$. Direct calculations lead to 
	$$
	\sum_{\nu} V(J_n, \phi_{\nu})^2 = \sum_{\nu} \frac{\omega^{2a + 1}}{(1 + \lambda\rho_{\nu})^2} \left\{\frac{1}{\sqrt{n}} \sum_{i = 1}^n [\zeta_i X_{i\nu} - \E(\zeta_i X_{i\nu})]\right\}^2,
	$$
	where $\zeta_i = \bK({-\eps_i(\theta_0)}/{h}) - \tau$ and $X_{i\nu} = \int_0^1 X_{i}(t)\phi_{\nu}(t)dt$ for $i = 1, \ldots, n$. As $|\zeta_i| \leq C$ holds almost surely for some positive constant $C > 0$, 
	\begin{align*}
		& \quad	\E\left\{\frac{1}{\sqrt{n}} \sum_{i = 1}^n [\zeta_i X_{i\nu} - \E(\zeta_i X_{i\nu})]\right\}^2 \\
		& = \Var(\zeta_iX_{i\nu}) \leq \E(\zeta_i^2X_{i\nu}^2)  \leq C^2 \E(X_{i\nu}^2) \\
		& \lesssim \E\|X_i\|_{L^2}^2 \cdot \int_0^1 \phi_{\nu}^2dt \\
		& \lesssim \nu^{2a}. 
	\end{align*}
	Furthermore, we have $\sum_{\nu} \frac{\omega^{2a + 1}\nu^{2a}}{(1 + \lambda\rho_{\nu})^2} < \infty$. Then, for any $\eps > 0$ and $\delta > 0$, there exists some $\nu_0$ such that $\sum_{\nu \geq \nu_0} \frac{\omega^{2a + 1}\nu^{2a}}{(1 + \lambda\rho_{\nu})^2} < \eps\delta$. 
	For any $\eps > 0$ and $\delta > 0$, we have by the Markov inequality, 
	$$
	\lim\sup_{n} P \left(\sum_{\nu \geq \nu_0} V(J_n, \phi_{\nu})^2 \geq \eps \right) < \lim\sup_{n}
	\frac{\E \left[\sum_{\nu \geq \nu_0} V(J_n, \phi_{\nu})^2\right] }{\eps} < \delta. 
	$$
	This shows (i) that $J_n$ is asymptotically finite-dimensional.
	
	Since $\eta(X_i) = \sum_{\nu} \frac{X_{i\nu}}{1 + \lambda\rho_{\nu}} \phi_{\nu}$, 
	$$
	V(J_n, \phi_{\nu}) = \frac{\omega^{a + 1/2}}{1 + \lambda\rho_{\nu}} \left\{\frac{1}{\sqrt{n}} \sum_{i = 1}^n [\zeta_i X_{i\nu} - \E(\zeta_i X_{i\nu})]\right\}. 
	$$
	Using the idea of proving Corollary \ref{cor-pointwise}, we can easily show that the Lindberg's condition is satisfied for the triangular array $\frac{1}{\sqrt{n}} \sum_{i = 1}^n [\zeta_i X_{i\nu} - \E(\zeta_i X_{i\nu})]$. Let $\sigma_{n\nu}^2 = \Var(\zeta_iX_{i\nu})$. Note that $\sigma_{n\nu}^2 \asymp 1$. Its proof is similar to the argument of $\Var(\zeta_i \inner{R_{X_i}}{\theta_1})$ in the proof of Corollary \ref{cor-pointwise}. Actually, as $n \rightarrow \infty$, $\sigma_n^2$ will converge to $\psi_{\nu}^2 := (\tau - \tau^2)\E[\int_0^1 X_i(t)\phi_{\nu}(t)dt]^2$. 
	Then $V(J_n, \phi_{\nu}) \convd N\left(0, \frac{\psi_{\nu}^2\omega^{2a + 1}}{(1 + \lambda\rho_{\nu})^2} \right)$. 
	
	By the Karhunen-Lo\`eve theorem, the Gaussian process $\ca G$ can be written as 
	$$
	\ca G(t) = \sum_{\nu} \frac{\omega^{a+1/2}\psi_{\nu}}{(1 + \lambda\rho_{\nu})} \xi_{\nu} \phi_{\nu}(t), 
	$$
	where $\xi_{\nu}$'s are i.i.d standard normal random variables. Hence $V(\ca G, \phi_{\nu})$ follows $N\left(0, \frac{\psi_{\nu}^2\omega^{2a + 1}}{(1 + \lambda\rho_{\nu})^2} \right)$. Namely, $V(J_n, \phi_{\nu})$ converges in distribution to $V(\ca G, \phi_{\nu})$. Then (ii) finite-dimensional convergence is verified. 
	
	This completes the proof.

	\subsection{Proof of Theorem \ref{th-quantCI}}
	
	We first treat the main term $\hat{Q}_{Y|X}(\tau | x_0) - Q_{Y | X}(\tau | x_0)$. Note that for a fixed $\tau \in \mU$, 
	\begin{align*}
		\hat{Q}_{Y|X}(\tau | x_0) - Q_{Y | X}(\tau | x_0) & = \hat{\alpha}_{h,\lambda} + \int_0^1 x_0(t)\hat{\beta}_{h,\lambda}(t, \tau)dt - \alpha_0(\tau) - \int_0^1 x_0(t)\beta(t,\tau) dt \\
		& = \inner{R_{x_0}}{\hat{\theta}_{h,\lambda} - {\theta}_0 }\\
		& = \inner{R_{x_0}}{\hat{\theta}_{h,\lambda} - {\theta}_0 + S_{n,h\lambda}(\theta_0)} - \inner{R_{x_0}}{S_{n,h\lambda}(\theta_0)}. 
	\end{align*}
	The first term is bounded by $\|R_{x_0}\| \cdot \|\hat{\theta}_{h,\lambda} - {\theta}_0 + S_{n,h\lambda}(\theta_0)\|$, which is of order $O_P(a_n\sigma_n(x_0))$ under the condition $\sigma_n(x_0) \asymp \|R_{x_0}\|$. We next derive the asymptotic distribution of the second term.

	Recall that $S_{n,h\lambda}(\theta_0) \Delta \theta  =    \frac{1}{n} \sum_{i = 1}^n\left[\bK\left(\frac{-\eps_i(\theta_0)}{h}\right) - \tau\right]\inner{R \lo {X_i}}{\Delta \theta}+ \lambda J(\beta_0, \Delta \beta)$ for any $\Delta\theta \in \ca H$. 	Let $\zeta_i = \bK({-\eps_i(\theta_0)}/{h}) - \tau, i = 1, \ldots, n$. Then
	\begin{align*}
		& \quad \inner{R_{x_0}}{S_{n,h\lambda}(\theta_0)} \\
		& = \frac{1}{n} \sum_{i = 1}^n \zeta_i ([\E\{B(X)\}]^{-1} + \inner{\eta(x_0)}{\eta(X_i)}_1) + \inner{P_{\lambda}\theta_0}{R_{x_0}} \\
		& = \frac{1}{n} \sum_{i = 1}^n \zeta_i Z_i + \int_0^1 x_0(t)(W_{\lambda}\beta_0)(t) dt,
	\end{align*}
	where $Z_i = [\E\{B(X)\}]^{-1} + \inner{\eta(x_0)}{\eta(X_i)}_1$ for $i = 1, \ldots, n$. To establish the CLT for the second term, we need to calculate $s_{2n}^2 = \Var\{\sum_{i = 1}^n [\bK({-\eps_i(\theta_0)}/{h}) - \tau] \inner{R_{x_0}}{R_{X_i}}\}$ first. 
	We employ the law of iterated expectations to find $\E \zeta_i^2Z_i^2$. As shown in the proof of Corollary \ref{cor-pointwise}, 
	$\E\{|\bK_h(-\eps(\theta_{0})) - \tau|^2 | X = x\}  = \tau - \tau^2 - 2hB(x) \int_0^{\infty} \bK(z) [1 - \bK(z)]dz + O(h^2)$ uniformly for $x$. By the Fourier expansions of $\eta(x_0)$ and $\eta(X_i)$, we have $\inner{\eta(x_0)}{\eta(X_i)}_1 = \sum_{\nu} \frac{x_{\nu}\hi 0 X_{i\nu}}{1 + \lambda\rho\lo{\nu}}$, where
	$x_{\nu}\hi 0 = \int_0^1 x_0(t)\phi_{\nu}(t)dt$ and $X_{i\nu} = \int_0^1 X_i(t)\phi_{\nu}(t)dt$ for $\nu \geq 1$. 
	Therefore, by Assumption \ref{ass-f} (c), 
	\begin{align*}
		\E \zeta_i^2 Z_i^2 & \leq C_1 \E \zeta_i^2 B(X_i) Z_i^2 \\
		& \asymp C_1 (\tau - \tau^2) \E\{B(X_i) [[\E\{B(X)\}]^{-2}  \\
		&~~~~~~~~~+ 2[\E\{B(X)\}]^{-1}\inner{\eta(x_0)}{\eta(X_i)}_1 + \inner{\eta(x_0)}{\eta(X_i)}_1^2]\} \\
		& = C_1 (\tau - \tau^2)\left([\E\{B(X)\}]^{-1} + \sum_{\nu} \frac{|x_{\nu}\hi 0|^2}{(1 + \lambda\rho\lo{\nu})^2}\right) \\
		& \asymp \sigma_n^2(x_0). 
	\end{align*}
	Next we find the $\E [\zeta_i \inner{\eta(x_0)}{\eta(X_i)}_1]$. It is actually $R^{(1)}_h(\theta_0; \tau)\theta_1$ if we take $\theta_1 = ([\E\{B(X)\}]^{-1} , \eta(x_0))$. 
	By the proof of (ii) of Lemma \ref{lemma1}, we have
	\begin{align*}
		|R^{(1)}_h(\theta_0; \tau)\theta_1 - R^{(1)}(\theta_0; \tau)\theta_1| & \leq Ch^{s+1} \E[|\inner{R \lo {X}}{ \theta_1}|]  \\
		& \leq  Ch^{s+1} \left[\E\left(\langle R_X, \theta_1 \rangle^4\right)\right]^{1/4} \\
		& \leq C h^{s+1} \|\theta_1\| = C h^{s+1} \|R_{x_0}\|
	\end{align*}
	by the Cauchy-Schwarz inequality and Lemma S.3 of \cite{shang2015b}. On the other hand, $\|R_{x_0}\| \asymp \sigma_n(x_0)$ by the condition. 
	Since $R^{(1)}(\theta_0; \tau)\theta_1 = 0$, $\{\E [\zeta_i \inner{\eta(x_0)}{\eta(X_i)}_1]\}^2 = O(h^{2s+2}) \sigma_n^2(x_0)$. 
	On the other hand, since $h = o(1)$, $\{\E [\zeta_i Z_i]\}^2= o(\E [\zeta_i^2Z_i^2])$. Hence we conclude that $s_{2n}^2 \asymp n\sigma_n^2(x_0)$. It implies that there exists a constant $\sigma_0^2 > 0$ that does not depend on $n$, such that $s_{2n}^2 \geq \sigma_0^2 n\sigma_n^2(x_0)$ for all $n \geq 1$.  
	
	We check the Lindeberg's condition to establish the central limit theorem for the triangular array $\sum_{i = 1}^n \{\zeta_i Z_i - R^{(1)}_h(\theta_0; \tau)\theta_1 \}$, where $\theta_1 = R_{x_0} = ([\E\{B(X)\}]^{-1}, \eta(x_0)) \in \ca H$. 
	Rearranging $\zeta_i Z_i - R^{(1)}_h(\theta_0; \tau)\theta_1$ leads to
	\begin{align*}
		\zeta_i Z_i - R^{(1)}_h(\theta_0; \tau)\theta_1 & = \left[\int_{-\infty}^{\inner{R \lo {X_i}}{ \theta_0}} k_h(u - Y_i) du - \tau\right] \inner{R \lo {X_i}}{ \theta_1} - \\
		& \quad\quad \int \left\{\int_{-\infty}^{\inner{R \lo {x}}{ \theta_0}} (\E[k_h(u - Y | x)]) du - \tau\right\} \inner{R \lo {x}}{ \theta_1} d P_X(x) .
	\end{align*}
	Obviously, there exists a positive constant $C$ that is independent of $n$ such that $|{\zeta}| \leq C$ a.s. As shown in the proof of Corollary \ref{cor-pointwise}, $|Z_i| = |\inner{R \lo {X}}{ \theta_1}| \leq  [\E\{B(X)\}]^{-1} + C_R^{1/2} \|\eta(x_0)\|_1 \|X\|_{L^2} \omega^{-(2a+1)/2}$ for some constant $C_R > 0$. Therefore, \\
	$|Z_i|^2 \lesssim (1 + C_R\omega^{-(2a + 1)} \|X_i\|_{L^2}^2) \|R_{x_0}\|^2$.
	
	Since $\log(\omega^{-1}) = O(\log n)$, we can choose $\ti{C}$ sufficiently large such that $\ti{C} \log n/ s  > C$ and $\omega^{-(2a + 1)}n^{-\ti{C}} = o(1)$. 
	Using the same idea of the proof of Corollary \ref{cor-pointwise}, we have for any $\epsilon > 0$, as $n$ approaching $\infty$, 
	\begin{align*}
		& \quad	ns_{2n}^{-2} \E \{[\zeta_i Z_i - R^{(1)}_h(\theta_0; \tau)\theta_1| ^2 \bbone(|\zeta_i Z_i - R^{(1)}_h(\theta_0; \tau)\theta_1|^2 \geq \eps^2s_{2n}^2)  \} \\
		& \lesssim \frac{\|R_{x_0}\|^2}{\sigma_n^2(x_0)} (\E\{(1 + C_R\omega^{-(2a + 1)} \|X_i\|_{L^2}^2)^2\})^{1/2} \cdot P(|\zeta_i|^2 Z_i^2 \geq \eps^2 s_{2n}^2)^{1/2} \\
		& \leq O(\omega^{-(2a+1)}) \left[P(s|\zeta_i| \geq \ti{C}\log n) + P\left(Z_i^2 \geq \frac{s^2\eps^2s_{2n}^2}{(\ti{C}\log n)^2}\right)\right]^{1/2} \\
		& \leq O(\omega^{-(2a+1)}) \left[P\left(s\|X\|_{L^2}^2 \geq s\sqrt{\frac{\omega^{2a+1}}{C_R} 
			\left(\frac{s^2\eps^2s_{2n}^2}{(\ti{C}\log n)^2 \|R_{x_0}\|^2}-1\right)}    \right)\right]^{1/2} \\
		& \leq O(\omega^{-(2a+1)}) \exp \left(-\frac{s}{2}\sqrt{\frac{\omega^{2a+1}}{C_R} 
			\left(\frac{s^2\eps^2s_{2n}^2}{(\ti{C}\log n)^2 \|R_{x_0}\|^2}-1\right)} \right) \\
		& = o(1),
	\end{align*}
	where the last equation follows by $s_{2n}^2 \asymp n \sigma_n^2(x_0) \asymp n \|R_{x_0}\|^2$, the choice of $\ti{C}$ and $n\omega^{2a+1} (\log n)^{-4} \rightarrow \infty$. 
	Thus the Lindberg's condition is fulfilled. 
	As a result, 
	$s_{2n}^{-1} \sum_{i = 1}^n \{\zeta_i Z_i - \E[\zeta_i Z_i]\} \convd N(0, 1)$. Moreover, $s_{2n}^{-1}n   \E[\zeta_i Z_i] = O(\sqrt{n} \sigma_n^{-1}(x_0) \cdot h^{s + 1} \|R_{x_0}\|) = O(\sqrt{n} h^{s + 1}) = o(1)$. Therefore, following the condition that $n a_n \sigma_n(x_0) /s_{2n} \asymp n^{1/2}a_n = o(1)$, we get that as $n \rightarrow \infty$, 
	\begin{align*}
		& \quad \frac{n}{s_{2n}} \left[\hat{Q}_{Y|X}(\tau | x_0) - Q_{Y | X}(\tau | x_0) + \int_0^1 x_0(t) (W_{\lambda}(\beta_0))(t)dt \right]  \\
		& = s_{2n}^{-1} \sum_{i = 1}^n \{\zeta_i Z_i - \E[\zeta_i Z_i]\} + ns_{2n}^{-1}  \E[\zeta_i Z_i]  + O_P(n a_n \sigma_n(x_0)/s_{2n}) \\
		& = s_{2n}^{-1} \sum_{i = 1}^n \{\zeta_i Z_i - \E[\zeta_i Z_i]\}  + o(1) + o_P(1) \convd N(0, 1). 
	\end{align*}
	
	Recall that if $\beta_0(\cdot, \tau) = \sum_{\nu} b_{\nu}\phi_{\nu}$,  
	$W_{\lambda}\beta_0 = \sum_{\nu} b_{\nu} \frac{\lambda\rho_{\nu}}{1 + \lambda\rho_{\nu}}\phi_{\nu}$. By Cauchy inequality the bias term satisfies that
	\begin{align*}
		& \quad	\bigg| \int_0^1 x_0(t) (W_{\lambda}(\beta_0))(t)dt \bigg| \\
		& = \bigg| \sum_{\nu} b_{\nu} \frac{\lambda\rho_{\nu}}{1 + \lambda\rho_{\nu}} \int_0^1 x_0(t) \phi_{\nu}(t) dt \bigg| \\
		& = \bigg| \sum_{\nu} b_{\nu} \frac{\lambda\rho_{\nu}}{1 + \lambda\rho_{\nu}} x_{\nu}^0 \bigg| \\
		& \leq \lambda (\sum_{\nu} b_{\nu}^2 \rho_{\nu}^2)^{1/2} \left(\sum_{\nu} \frac{|x_{\nu}^0|^2}{(1 + \lambda\rho_{\nu})^2}\right)^{1/2}. 
	\end{align*}
	On the other hand, since $s_{2n}^2 \geq \sigma_0^2 n \sigma_n^2(x_0) \geq \sigma_0^2 n \sum_{\nu} \frac{|x_{\nu}^0|^2}{(1 + \lambda\rho_{\nu})^2}$ and $\sum_{\nu} b_{\nu}^2 \rho_{\nu}^2 < \infty$ and $n\omega^{4l} = o(1)$, we have
	\begin{align*}
		& \quad	(n/s_{2n}) \bigg| \int_0^1 x_0(t) (W_{\lambda}(\beta_0))(t)dt \bigg| \\ 
		& \leq \frac{n}{\sqrt{\sigma_0^2 n \sum_{\nu} \frac{|x_{\nu}^0|^2}{(1 + \lambda\rho_{\nu})^2}}} \cdot \lambda (\sum_{\nu} b_{\nu}^2 \rho_{\nu}^2)^{1/2} \cdot \left(\sum_{\nu} \frac{|x_{\nu}^0|^2}{(1 + \lambda\rho_{\nu})^2}\right)^{1/2} \\
		& = O(n^{1/2}\lambda) = o(1).
	\end{align*}
	Hence, 
	$$
	(n/s_{2n}) \left[\hat{\alpha}_{h,\lambda} + \int_0^1 x_0(t) \hat{\beta}_{h,\lambda}(t, \tau)dt - \alpha_0(\tau) - \int_0^1 x_0(t)\beta(t,\tau) dt \right] \convd N(0, 1)
	$$
	as $n \rightarrow \infty$.

	\bibliography{sfqr}{}
	\bibliographystyle{natbib}

\end{document}